\documentclass{IEEEtran}
\usepackage{cite}
\usepackage{amsmath,amssymb,amsfonts,amsthm}
\usepackage{algorithmic}
\usepackage{graphicx}
\usepackage{textcomp}
\usepackage{bm}
\usepackage[lined,ruled,vlined,boxed]{algorithm2e}
\usepackage{booktabs}
\usepackage{color}
\usepackage{nomencl}
\usepackage{multirow}
\usepackage{threeparttable}
\usepackage{enumitem} 
\usepackage{nomencl}
\usepackage{subfig} 
\makenomenclature

\newtheorem{lemma}{Lemma}
\newtheorem{theorem}{Theorem}
\newtheorem{assumption}{Assumption}
\newtheorem{remark}{Remark}
\newtheorem{proposition}{Proposition}
\newtheorem{definition}{Definition}

\def\BibTeX{{\rm B\kern-.05em{\sc i\kern-.025em b}\kern-.08em
		T\kern-.1667em\lower.7ex\hbox{E}\kern-.125emX}}
\usepackage{balance}

\allowdisplaybreaks[3]

\begin{document}
\title{Perturbed Proximal Gradient ADMM for Nonconvex Composite Optimization}
\author{Yuan Zhou, Xinli Shi, \IEEEmembership{Senior Member, IEEE}, Luyao Guo, Jinde Cao, \IEEEmembership{Fellow, IEEE} and Mahmoud Abdel-Aty
		

		

}
\maketitle
		
	\begin{abstract}
 This paper proposes a Perturbed Proximal Gradient ADMM (PPG-ADMM) framework for general nonconvex composite optimization problems.
 The core of PPG-ADMM is a perturbation mechanism that enhances convergence stability and, more importantly, relaxes two restrictive assumptions common in existing ADMM-based methods: smoothness of the final objective block and range inclusion conditions on the constraint matrices.
Convergence is established through a technically constructed Lyapunov function, which guarantees sufficient descent and has a well-defined lower bound.
With appropriate parameters, PPG-ADMM converges to an $\epsilon$-approximate stationary point with a sublinear convergence rate of $\mathcal{O}(1/\sqrt{K})$.
Furthermore, with a suitably decaying perturbation parameter, the algorithm can attain an $\epsilon$-stationary point, ensuring stronger optimality guarantees.
The framework is extended to fully decentralized settings, and its effectiveness is validated on a distributed nonconvex composite optimization task.
	\end{abstract}
	

	\section{Introduction}\label{Sec1}

The Alternating Direction Method of Multipliers (ADMM) \cite{boyd2011distributed} has been widely successful in convex optimization. 
However, many contemporary applications, such as smart grid management \cite{zhou2023proximal}, matrix separation \cite{shen2014augmented}, and sparse principal component analysis \cite{hajinezhad2015nonconvex}, are formulated as Nonconvex Composite Optimization Problems (NCOPs). 
These problems typically involve coupled nonconvex and nonsmooth terms, rendering standard convex convergence analysis inapplicable. 
This creates an urgent need for ADMM frameworks specifically tailored to the theoretical challenges of NCOPs.
	
	We consider the following linearly constrained problem:
\begin{equation}\label{ClassicalADMMProblem}
\min_{\mathbf{x},\mathbf{z}}  \  F(\mathbf{x})+H(\mathbf{z}), \quad s.t. \  \mathbf{A}\mathbf{x}+\mathbf{B}\mathbf{z}=\mathbf{c},
	\end{equation}
where $\mathbf{x}\in\mathbb{R}^n$, $\mathbf{z}\in\mathbb{R}^{m}$, $\mathbf{A}\in\mathbb{R}^{p\times n}$,  $\mathbf{B}\in\mathbb{R}^{p\times m}$, $\mathbf{c}\in\mathbb{R}^{p}$, $F=F^0+F^1$ and $H=H^0+H^1$ admit a decomposition into smooth nonconvex and nonsmooth weakly convex components \cite{barber2024convergence,han2018linear}.
Establishing convergence for such general NCOPs is challenging. 
Specifically, standard analysis typically relies on constructing a Lyapunov function based on the Lagrangian.
However, to ensure sufficient descent, existing works typically impose two restrictive but necessary conditions
\cite{sun2019two,yang2022proximal,yang2022survey,wang2019global}:
	\begin{enumerate}[]
		\item[C1:] The matrices satisfy $\mathrm{Im}(\mathbf{A})\cup\{\mathbf{c}\}\subseteq \mathrm{Im}(\mathbf{B})$.
		\item[C2:] The objective term associated with the last updated primal variable is smooth.
	\end{enumerate}
    \begin{table*}[thb]
 \renewcommand\arraystretch{1.15}
		\centering
  \begin{threeparttable}
		\caption{{ Some Existing ADMMs for Nonconvex Composite Optimization.}}\label{table1}
		\setlength{\tabcolsep}{8pt}
		\begin{tabular}{c|c|c|c|c}
			\toprule
			{Reference} & {Objective Function\tnote{1}}  & {Constraint\tnote{1}} &Smoothness  & {Assumptions of Matrices in Constraints}\\
   \midrule
   \cite{wang2019distributed,chang2016asynchronous} & \multirow{6}{*}{$F(x)+H(z)$} &$x-z=0$ &\multirow{6}{*}{$H$} & —— \\
\cite{li2015global,bian2021stochastic,boct2020proximal} &  & $x-Bz=0$ & &$B$ is full row rank. \\
\cite{guo2017convergence} & &  $Ax+z=c$ & &$A$ is full column rank. \\
{\cite{huang2018mini,huang2019faster}} & &  {$Ax+Bz=c$} & &{$B$ is full row or column rank}. \\
{\cite{zeng2024unified,zeng2024accelerated}} & &  {$Ax+Bz=c$} & &{$\mathrm{Im}(A)\cup \{c\}\subseteq\mathrm{Im}(B)$}. \\
\cite{yashtini2021multi} & &  $Ax+Bz=c$ & &$\mathrm{Im}(A)\subset\mathrm{Im}(B)$ and $c\in \mathrm{Im}({B})$. \\
\hline
{\cite{hong2016convergence,hong2017distributed}\tnote{2}} &\multirow{3}{*}{$F(x)+G(x)+H(z)$}  & $x-z=0$ & {$H$} & —— \\
{\cite{hong2016convergence}\tnote{2}}& & $Ax-z=0$ & {$H$} & $A$ is full column rank. \\
\cite{yashtini2022convergence} & &  $Ax+Bz=c$ &$H$, $G$ & $B$ is full rank, $\mathrm{Im}({A})\subseteq \mathrm{Im}({B})$, $c\in \mathrm{Im}({B})$. \\
\hline
\cite{liu2019linearized} & $F(x)+G(x,z)+H(z)$ &$Ax+Bz=0$ &$H$, $G$ & $B$ is full column rank, $\mathrm{Im}(A)\subset\mathrm{Im}(B)$. \\
\hline
\cite{hajinezhad2018alternating} &\multirow{5}{*}{$F(x)+G(y)+H(z)$} &  $xy-z=0$ & \multirow{5}{*}{$H$} & —— \\
\cite{zhang2021extended} & & $A_1x+A_2y+z=c$ & & $A_2$ is full column rank. \\
\cite{yang2017alternating} & & $A_1x+A_2y-z=0$ & & $A_1$ and $A_2$ are full column rank. \\
\cite{wang2018convergence} &  &$A_1x+A_2y+Bz=0$ & &$B$ is full row rank. \\
\cite{melo2017iteration} &  &$A_1x+A_2y+Bz=c$ & &$B\neq \mathbf{O}$, $\mathrm{Im}(A_1)\cup\mathrm{Im}(A_2)\cup \{c\}\subset\mathrm{Im}(B)$. \\
\hline
\cite{bot2019proximal} &$F(x)+G(y)+H(y,z)$ & $x-Bz=0$ &$H$ &$B$ is full row rank. \\
			\hline
   This paper &$F^0(x)+F^1(x)+H^0(z)+H^1(z)$ & $Ax+Bz=c$ &$F^0, H^0$ & The feasible set is nonempty. \\
   \bottomrule
		\end{tabular}
  \begin{tablenotes} 
        \item[1] In this table, $z$ corresponds to the final updated primal variable in these algorithms, while $B$ corresponds to the matrix $\mathbf{B}$ mentioned in C1.
        To achieve a more unified representation of optimization problems, separable objective functions are expressed in a compact form, and some symbols used in certain literature have been substituted.
        \item[2] Set constraints of $x$ are replaced with indicator functions $G$ here.
        \vspace{-10pt}
     \end{tablenotes} 
     \end{threeparttable}
	\end{table*}
    Specifically, C1 ensures the feasibility of the linear constraints, while both conditions are indispensable for bounding the dual update difference (which typically appears as a non-negative term in the Lyapunov difference) using the primal residual to guarantee Lyapunov descent.
    As summarized in TABLE \ref{table1}, representative ADMM-based methods \cite{wang2019distributed,li2015global,bian2021stochastic,boct2020proximal,guo2017convergence,yashtini2021multi,hajinezhad2018alternating,yashtini2022convergence,zhang2021extended,yang2017alternating,wang2018convergence,bot2019proximal,liu2019linearized,hong2016convergence,melo2017iteration,huang2018mini,huang2019faster,zeng2024unified,zeng2024accelerated,chang2016asynchronous,hong2017distributed} universally require both conditions.
    However, these conditions are often violated in practice.
    For instance, in Robust Principal Component Analysis (RPCA) \cite{wen2019robust,feng2024online}, both objective terms are typically nonsmooth (violating C2), and constraint matrices may be rank-deficient (potentially violating C1).

    Motivated by these limitations, we propose a Perturbed Proximal Gradient ADMM (PPG-ADMM) framework. 
    Our algorithm integrates a perturbation mechanism into the dual update and employs linearization techniques combined with proximal operators for the primal subproblems. 
    The main contributions of this paper are summarized as follows:

    \textbf{Broader Applicability via Perturbation}. 
    The proposed perturbation mechanism allows PPG-ADMM to explicitly construct a valid Lyapunov function, circumventing the restrictive conditions C1 and C2. 
    This enables the algorithm to handle a wider class of NCOPs, significantly broadening applicability.
Furthermore, this mechanism strikes a balance between optimality and constraint satisfaction, which reduces the sensitivity of dual updates and mitigates oscillations, thereby enhancing convergence stability even in ill-conditioned problems.

    \textbf{Rigorous Convergence Analysis}. 
    We establish that PPG-ADMM achieves a convergence rate of $\mathcal{O}(1/\sqrt{K})$ and an iteration complexity of $\mathcal{O}(1/\epsilon^2)$ for reaching an $\epsilon$-AKKT point, matching the lower bound for first-order methods in nonconvex optimization.
By setting the perturbation parameter $\beta=\mathcal{O}(1/K)$, PPG-ADMM further attains an $\epsilon$-KKT point.
Crucially, these optimality guarantees are established for general NCOPs under mild conditions, without relying on restrictive compact-set assumptions.

\textbf{Decentralized Extension}. 
We further demonstrate the versatility of the PPG-ADMM framework by extending it to decentralized settings. 
We construct fully decentralized algorithms that operate without a central coordinator, thereby broadening the potential utility in networked systems.
Numerical experiment validates the effectiveness and robustness of the proposed algorithm. 

\subsection*{Related Work}
 We summarize a range of representative works \cite{wang2019distributed,li2015global,bian2021stochastic,boct2020proximal,guo2017convergence,yashtini2021multi,hajinezhad2018alternating,yashtini2022convergence,zhang2021extended,yang2017alternating,wang2018convergence,bot2019proximal,liu2019linearized,hong2016convergence,melo2017iteration,huang2018mini,huang2019faster,zeng2024unified,zeng2024accelerated,chang2016asynchronous,hong2017distributed} that utilize ADMM to solve NCOPs in TABLE \ref{table1}.
 This table highlights the assumptions corresponding to C1 and C2, illustrating that despite differences in problem formulations and algorithms, these works all impose both conditions.
 Several works \cite{wang2019distributed,chang2016asynchronous,li2015global,bian2021stochastic,boct2020proximal,guo2017convergence,huang2018mini,huang2019faster,yashtini2022convergence,liu2019linearized,zhang2021extended,bot2019proximal,wang2018convergence,yang2017alternating,hong2017distributed,hajinezhad2018alternating} make explicit assumptions on the structure of the constraint matrices, such as requiring certain matrices to be full-rank or identity matrices. 
 These are special cases of C1 and significantly restrict applicability, particularly in decentralized optimization where mixing matrices are typically rank-deficient.
To establish C2, all these works \cite{wang2019distributed,li2015global,bian2021stochastic,boct2020proximal,guo2017convergence,yashtini2021multi,hajinezhad2018alternating,yashtini2022convergence,zhang2021extended,yang2017alternating,wang2018convergence,bot2019proximal,liu2019linearized,hong2016convergence,melo2017iteration,huang2018mini,huang2019faster,zeng2024unified,zeng2024accelerated,chang2016asynchronous,hong2017distributed} require the function associated with the final updated primal variable, typically denoted as $H$, to be smooth, though they impose no strict conditions on its convexity or the smoothness of other terms. 
For example, $H$ is strongly convex in \cite{melo2017iteration}, whereas in \cite{hong2016convergence,chang2016asynchronous,wang2019distributed,hong2017distributed,huang2019faster,zeng2024accelerated,gao2024non,li2018simple,jiang2022distributed,wang2021distributed,chen2021communication,xin2021stochastic,mancino2023proximal,xiao2023one,liu2023proximal}, $F$ or $G$ can be convex and nonsmooth, with some of the nonsmooth terms are required to be weakly convex in \cite{chen2021distributed,yan2023compressed,bohm2021variable,liu2023proximal}.
When the proximal operators of these nonsmooth terms are easily computable, closed-form solutions can be obtained for the subproblems.
Otherwise, the lack of such favorable properties may lead to nonconvex and nonsmooth subproblems, making them difficult to solve.
When satisfying C1 and C2 simultaneously is infeasible, a two-level framework combining ADMM and the augmented Lagrangian method (ALM) is proposed in \cite{sun2019two}, where the original problem is relaxed and solved via an inner ADMM loop, with an outer ALM gradually reducing the relaxation parameter.
However, the complexity and computational burden of this nested-loop structure restrict its practical applicability.
In summary, designing a simple, efficient, and generally applicable ADMM for NCOPs without relying on C1 and C2 remains a valuable and open research topic.

The works in \cite{wang2019distributed,li2015global,bian2021stochastic,boct2020proximal,guo2017convergence,yashtini2021multi,hajinezhad2018alternating,yashtini2022convergence,zhang2021extended,yang2017alternating,wang2018convergence,bot2019proximal,liu2019linearized,hong2016convergence,melo2017iteration,huang2018mini,huang2019faster,zeng2024unified,zeng2024accelerated,chang2016asynchronous,hong2017distributed} exhibit several common strategies in the construction and analysis of ADMM for solving NCOPs.
In problems involving separable objective functions, Gauss-Seidel update schemes are commonly adopted  \cite{yashtini2021multi,wang2018convergence}, while other works \cite{liu2019linearized,wang2019distributed,hong2016convergence} allow multi-block parallel updates, significantly accelerating computation.
To reduce per-iteration complexity, \cite{liu2019linearized,yashtini2022convergence,boct2020proximal,bot2019proximal,huang2018mini,huang2019faster,zeng2024unified,zeng2024accelerated} utilize linear approximations to substitute smooth terms in the objective functions of subproblems, while \cite{li2015global,wang2018convergence,melo2017iteration,zhu2024first} further integrate Bregman distances.
For cases where certain terms exhibit convex or concave properties, \cite{hajinezhad2018alternating} constructs distinct surrogate functions tailored to these structures.
Beyond these structural modifications, many works introduce proximal terms into the ADMM subproblems to further improve numerical stability and avoid costly matrix inverse computations \cite{bian2021stochastic,yashtini2021multi,yashtini2022convergence,boct2020proximal,liu2019linearized,zhang2021extended,bot2019proximal,zeng2024accelerated,zeng2024unified,huang2018mini,huang2019faster,hajinezhad2018alternating}.
Strong convexity is typically enforced either through sufficiently large penalty parameters in the augmented Lagrangian \cite{hong2016convergence,chang2016asynchronous,wang2019distributed,hong2017distributed} or via the integration of Bregman distances \cite{li2015global,melo2017iteration}.
These strong convexity properties \cite{wang2019distributed,hong2016convergence,li2015global,melo2017iteration}, together with smoothness assumptions \cite{guo2017convergence,zhang2021extended}, facilitate the design of monotonically decreasing Lyapunov functions, laying the foundation for theoretical convergence guarantees \cite{wang2019distributed,hong2016convergence}.
In large-scale machine learning applications, stochastic ADMM variants combined with variance reduction techniques, such as SAGA, SVRG, and SPIDER, are explored in \cite{huang2018mini,huang2019faster,zeng2024unified,zeng2024accelerated}.
They achieves convergence to an $\epsilon$-stationary point (under Definition \ref{def:AKKT} provided in this paper) at a sublinear rate of $\mathcal{O}(1/\sqrt{K})$ for NCOPs, consistent with the theoretical result in \cite{melo2017iteration}.
Finally, several works in \cite{yashtini2021multi,yashtini2022convergence,yang2017alternating,wang2018convergence,bot2019proximal,zhu2024first,zeng2024accelerated,zeng2024unified} strengthen their theoretical contributions by leveraging the Kurdyka-Łojasiewicz (KL) condition, which provides refined convergence analysis under mild conditions.

In addition to ADMM, various other algorithms have been developed for solving NCOPs.
A projected subgradient method for optimization problems with weakly convex objective functions and set constraints is introduced in \cite{chen2021distributed}.
The Proximal Gradient Method (PGM) \cite{gao2024non,li2018simple,jiang2022distributed} is another commonly used approach.
The algorithm in \cite{gao2024non} combines PGM with Polyak's momentum, allowing for inexact solutions in the proximal step, while \cite{li2018simple} incorporates variance reduction techniques into stochastic PGM, establishing linear convergence under the general Polyak-Łojasiewicz (PL) condition. 
To handle nonsmooth and weakly convex terms involving linear operators, \cite{liu2023proximal} extends the variable smoothing technique of \cite{bohm2021variable} and proposes the Proximal Variable Smoothing Gradient (ProxVSG) method. 
The methods in \cite{liu2023proximal,bohm2021variable,yan2023compressed} approximate nonsmooth terms using smooth Moreau envelopes, which facilitate optimization but alter the original objective function and inevitably lead to approximate solutions.
Several primal-dual methods have also been proposed based on the Lagrangian framework, such as those in \cite{wang2021distributed,chen2021communication,zhu2024first}, offering flexible update schemes. 
For instance, \cite{zhu2024first} considers a setting where both C1 and C2 hold, and adopts a problem formulation closely related to that in \cite{liu2019linearized}.
However, many of these works \cite{chen2021distributed,gao2024non,li2018simple,jiang2022distributed,wang2021distributed,chen2021communication,yan2023compressed} consider relatively simple problem structures with only a single primal variable, which limits their applicability to more complex settings.
In addition, Proximal Alternating Linearized Minimization (PALM) methods \cite{wang2023generalized,pock2016inertial} have been applied to NCOPs, but they are restricted to unconstrained problems.

Most of the aforementioned algorithms for NCOPs are designed for centralized settings, whereas several works \cite{wang2019distributed,hong2016convergence,hong2017distributed,chang2016asynchronous,yan2023compressed,wang2021distributed,chen2021communication,jiang2022distributed,xin2021stochastic,mancino2023proximal,xiao2023one} explore distributed algorithms.
Among them, the distributed ADMM variants proposed in \cite{wang2019distributed,hong2016convergence,hong2017distributed,chang2016asynchronous} allow asynchronous updates across nodes using cyclic or randomized rules, removing the need for global synchronization.
However, these approaches typically rely on a star topology, where a central node coordinates communication and computation.
This introduces the risk of a single point of failure and imposes significant communication and computational burdens on the central node \cite{boyd2011distributed}.
To overcome these limitations, decentralized optimization methods have been proposed.
Literature \cite{jiang2022distributed} develops a decentralized PGM algorithm for time-varying networks by introducing dynamic mixing matrices.
The works \cite{wang2021distributed,chen2021communication} initially consider NCOPs with only $x$-related terms, but incorporate additional primal variables through proximal regularization to ensure strong convexity and facilitate algorithmic design and convergence analysis.
Further, decentralized stochastic Proximal Gradient Tracking (PGT) methods have been investigated in \cite{xin2021stochastic,yan2023compressed,mancino2023proximal,xiao2023one}, aiming to mitigate the adverse effects of data heterogeneity on convergence. 
Notably, some of these methods \cite{xin2021stochastic,mancino2023proximal,xiao2023one} demonstrate \textit{linear speedup}, where increasing the number of agents improves convergence rates. 
Additionally, \cite{yan2023compressed,chen2021communication} incorporate compression techniques to reduce communication overhead within the network.
Despite these advances, existing decentralized algorithms \cite{wang2021distributed,chen2021communication,jiang2022distributed,xin2021stochastic,mancino2023proximal,xiao2023one} are primarily restricted to consensus-type constraints and convex nonsmooth terms. 
In particular, decentralized ADMM-based methods applicable to general NCOPs remain limited in the current literature.

 \textbf{Notations:} 
In this paper, the $n$-dimensional vector spaces and the $n\times m$ matrix spaces are symbolized as $\mathbb{R}^n$ and $\mathbb{R}^{n\times m}$, respectively.
 The $\ell_1$ norm and $\ell_2$ norm are represented by $\Vert \cdot \Vert_1$ and $\Vert \cdot \Vert$, respectively. 
 The induced norm associated with a positive semidefinite matrix $\mathbf{M}$ is denoted as ${\Vert \cdot \Vert}_\mathbf{M}=\sqrt{{\langle\cdot,\mathbf{M}\cdot\rangle}}$. 
 The image/range space of matrix $\mathbf{M}$ is represented as $\mathrm{Im}(\mathbf{M})$.
 The inverse matrix and transpose of $\mathbf{M}$ are represented by $\mathbf{M}^{-1}$ and $\mathbf{M}^{\top}$, respectively.
	Additionally, $\mathbf{M}\succ(\succeq) \mathbf{N}$ indicates that $\mathbf{M}-\mathbf{N}$ is a positive (semi)-definite matrix.
	Furthermore, $\mathbf{0}$ represents vectors consisting of all zeros, and $\mathbf{O}$ and ${\mathbf{I}}$ are zero matrix and identity matrix, respectively.
 The differential operator and subdifferential operator are denoted by $\nabla $ and $\partial$.
 Define $\mathrm{dist}(\mathbf{x},\mathcal{C})=\inf_{\mathbf{y}\in\mathcal{C}}\|\mathbf{x}-\mathbf{y}\|$ for a set $\mathcal{C}\in\mathbb{R}^n$.

\section{Preliminaries}
 We provide several definitions and lemmas here which will be utilized in the subsequent analysis.
 
\begin{definition}\label{DefinitionWeaklyconvex} 
		(Weakly Convex): The function $g:\mathbb{R}^n\rightarrow\mathbb{R}\cup\{+\infty\}$ is weakly convex if $g(\cdot)+\frac{\gamma}{2}{\|\cdot\|}^2$ is convex with $\gamma\geqslant0$.
	\end{definition}
 Some weakly convex functions can be used as regularization terms in various optimization problems \cite{zhu2024first,hajinezhad2016nestt}. 
 For instance, the minimax concave penalty (MCP) and the smoothly clipped absolute deviation (SCAD) penalty are often used as alternatives to the $\ell_1$ norm in machine learning \cite{bohm2021variable} and RPCA \cite{wen2019robust}, as they reduce bias and enhance solution stability.
 The MCP, which is $1/\theta$-weakly convex, is defined as
  \begin{align*}
\mathcal{M}_{\eta,\theta}(x)\!=\! \begin{cases}
    \eta|x|-{x^2}/{(2\theta)}, & |x|\leqslant\theta\eta, \\
    {\theta \eta^2}/{2}, & \textrm{otherwise.}
    \end{cases} 
 \end{align*}
 {The $1/(\xi-1)$-weakly convex SCAD penalty is given by
\begin{align*}
    \mathcal{S}_{\eta,\xi}(x)\!=\!\begin{cases}
\eta |x|, & |x| \leqslant \eta, \\
\frac{2\xi\eta |x|-x^2-\eta^2 }{ 2(\xi-1) }, & \eta < |x| \leqslant \xi\eta, \\
{(\xi+1)\eta^2}/{2}, & |x| > \xi\eta,
\end{cases}
\end{align*}
where $\theta>0$, $\xi>2$ and $\eta>0$ are tuning parameters that control the shape of the regularization.} 

	\begin{definition}\label{DefinitionProx} 
		(Proximal Operator): Consider a $\gamma$-weakly convex function $g:\mathbb{R}^n\rightarrow\mathbb{R}\cup\{+\infty\}$, for $x, y\in\mathbb{R}^n$ and $\tau>\gamma$, its proximal operator is defined as
  \begin{equation*}
   \textbf{prox}_{g}^{\tau}(y)=\arg \min_x\{g(x)+\frac{\tau}{2}{\Vert x-y \Vert}^2\}.
  \end{equation*}
	\end{definition}
 The function to be minimized above is strongly convex, ensuring that the corresponding proximal operator is well-defined and unique \cite{liu2023proximal,bohm2021variable}.
 The proximal operators of several commonly used nondifferentiable functions are well established.
For instance, the proximal operator of $\ell_1$-norm corresponds to the soft-thresholding operator, while the proximal operators of MCP and SCAD penalty are given by
  \begin{align*}
    \textbf{prox}_{\mathcal{M}_{\eta,\theta}}^{\tau}(x)\!=\! \begin{cases}
     0, & |x|<\eta/\tau,\\
    \frac{\tau\theta x-\text{sign}(x)\cdot\theta\eta }{\tau\theta-1}, & \eta/\tau\leqslant|x|\leqslant\theta\eta, \\
    x, & |x|>\theta\eta,
    \end{cases} 
 \end{align*}
 {\begin{align*}
     \textbf{prox}_{\mathcal{S}_{\eta,\xi}}^{\tau}(x)\!=\!
\begin{cases}
\text{sign}(x)\!\cdot\! (|x|\!-\!\eta/\tau)_+, & |x| \leqslant (1+1/\tau)\eta, \\
\frac{\tau(\xi-1)x-\text{sign}(x)\cdot\xi\eta}{\tau\xi-\tau-1}, & (1\!+\!1/\tau)\eta \!< \!|x|\!\leqslant\! \xi\eta, \\
x, & |x| > \xi\eta,
\end{cases}
 \end{align*}
where $(|x|-\eta/\tau)_+=\max\{|x|-\eta/\tau,0\}$ is the ReLU function.}

 \begin{definition}\label{DefinitionSub} 
	(Subdifferential): For a proper and closed function $g:\mathbb{R}^n\rightarrow\mathbb{R}\cup\{+\infty\}$. 
	The Fr\'{e}chet subdifferential $\hat{\partial} g$ and the limiting subdifferential $\partial g$ is denoted respectively as
	$$
	\hat{\partial} g(x)=\{v\in\mathbb{R}^n: \lim_{y\rightarrow x}\inf_{y\neq x} \frac{g(y)-g(x)-\langle v,y-x\rangle}{\|y-x\|}\geqslant0\},
	$$  
 	$$
	\partial g(x)=\{v\in\mathbb{R}^n: \exists \  x^k\stackrel{g}{\rightarrow} x, v^k\rightarrow v \ \text{with} \ v^k\in\hat{\partial} g(x^k)\},
	$$  
 where $x^k\stackrel{g}{\rightarrow} x$ means $x^k{\rightarrow} x$ and $g(x^k){\rightarrow} g(x)$.
\end{definition}
The Fr\'{e}chet subdifferential set is closed and convex, whereas the limiting subdifferential set is closed \cite{rockafellar2009variational}.
Furthermore, we have $\hat{\partial} g(x) \subseteq{\partial} g(x)$, and if $g$ is convex, then 
 	$$
	\hat{\partial} g(x) ={\partial} g(x)=\{v\in\mathbb{R}^n: g(y)-g(x)\geqslant\langle v, y-x \rangle\}.
	$$  
The subdifferential of a convex function is monotone \cite{ryu2022large}, i.e.,
 $$
 \langle v_x-v_y, x-y \rangle\geqslant 0, v_x\in\partial g(x), v_y\in\partial g(y).
 $$
 Further, if $g$ is $\gamma$-weakly convex, it further has
 \begin{equation}\label{weakConvexP}
\langle v_x-v_y, x-y \rangle\geqslant -\gamma{\|x-y\|}^2, v_x\in\partial g(x), v_y\in\partial g(y).
 \end{equation}
 It can be directly deduced by the convexity of $g(\cdot)+\frac{\gamma}{2}{\|\cdot\|}^2$.
 \begin{definition}\label{LemmaStronglyConvex}
		(Strong Convexity): The function $g:\mathbb{R}^n\rightarrow\mathbb{R}$ is called $\mu$-strongly convex with modulus $\mu>0$, for any $x, y\in\mathbb{R}^n$, it holds that
  \begin{subequations}
  \begin{align}
            &\qquad\langle v_x-v_y, x-y \rangle\geqslant \mu{\Vert x-y\Vert}^2, v_x\in\partial g(x), \label{LemmaStronglyConvexP1}\\
      &g(x)-g(y)\geqslant \langle v_y, x-y\rangle+\frac{\mu}{2}{\Vert x-y\Vert}^2, v_y\in\partial g(y).\label{LemmaStronglyConvexP2}
  \end{align}
  \end{subequations}
	\end{definition}
	
\begin{lemma}\cite[Theorem 5.17]{beck2017first}\label{LemmaStronglyConvex2} 
The function $g$ is $\mu$-strongly convex if and only if $g(\cdot)-\frac{\mu}{2}{\|\cdot\|}^2$ is convex.
	\end{lemma}
We can further deduce from Lemma \ref{LemmaStronglyConvex2} that if $f$ is $\mu$-strongly convex and $g$ is $\gamma$-weakly convex with $\mu>\gamma$, then $f(\cdot)+g(\cdot)-\frac{\mu-\gamma}{2}{\|\cdot\|}^2$ will be convex, which implies function $f+g$ is $(\mu-\gamma)$-strongly convex.
 
 \begin{definition}\label{DefinitionCoercive} 
		(Coercive Function): The function $g:\mathbb{R}^n\rightarrow\mathbb{R}\cup\{+\infty\}$ is coercive if it satisfies
		$
		\lim\limits_{\Vert x\Vert\rightarrow+\infty}\mathop{g(x)}=\infty.
		$
	\end{definition}
 \begin{definition}\label{DefinitionLip} 
	(Smoothness): The function $g$ is differentiable and has gradient Lipschitz continuous with modulus $L$, i.e.,
$$
\Vert\nabla  g(x)-\nabla  g(y)\Vert\leqslant L\Vert x-y\Vert, \ \forall x, y\in\mathbb{R}^n.
$$
Then, $g$ is called Lipschitz differentiable, or $L$-smooth.
\end{definition}
\begin{lemma}\cite{liu2019linearized}\label{LemmaDifferentiable} 
		Suppose that $f(x)$ is differentiable and $g(x)$ is possibly nondifferentiable. 
		Assume that there exists $v_1\in\partial g(x)$, we have
		$v_2=\nabla  f(x)+v_1\in\partial [f(x)+g(x)]$.
	\end{lemma}

	\section{Algorithm Development}
Prior to introducing PPG-ADMM framework for NCOP \eqref{ClassicalADMMProblem}, it is essential to introduce the assumptions made in this paper.
Similar assumptions  are prevalent in related literature \cite{li2015global,bian2021stochastic,boct2020proximal,guo2017convergence,yashtini2021multi,hajinezhad2018alternating,yashtini2022convergence,zhang2021extended,yang2017alternating,wang2018convergence,bot2019proximal,liu2019linearized,hong2016convergence,yang2022proximal,sun2019two,hajinezhad2019perturbed,wang2019global,wang2019distributed,melo2017iteration,huang2018mini,huang2019faster,zeng2024accelerated,zeng2024unified,chang2016asynchronous,hong2017distributed}, albeit with slight variations across works.
 \begin{assumption}[Feasibility] \label{AssumptionMatrix}
		Problem \eqref{ClassicalADMMProblem} is feasible and its set of stationary points is nonempty.
	\end{assumption}
\begin{assumption}\label{AssumptionFunction}
		The objective function in \eqref{ClassicalADMMProblem} satisfies:
		\begin{enumerate}[label=(\roman*)]
			\item  $F$ and $H$ are lower bounded and coercive;\label{AssumptionFunctioni}
			\item  $F^0$ and $H^0$ are $L_F$- and $L_H$-smooth, respectively;\label{AssumptionFunctionii}
			\item  $F^1$ and $H^1$ are $\gamma_F$- and $\gamma_H$-weakly convex, respectively.
		\end{enumerate}
	\end{assumption}

    The function $g:\mathbb{R}^n\rightarrow\mathbb{R}\cup\{+\infty\}$ is weakly convex if $g(\cdot)+\frac{\gamma}{2}{\|\cdot\|}^2$ is convex with $\gamma\geqslant0$.
 Some weakly convex functions can be used as regularization terms in various optimization problems. 
 For instance, the minimax concave penalty (MCP) and the smoothly clipped absolute deviation (SCAD) penalty are often used as alternatives to the $\ell_1$ norm in machine learning \cite{bohm2021variable} and RPCA \cite{wen2019robust}, as they reduce bias and enhance solution stability.
 
  

 The algorithm framework is built upon the following perturbed augmented Lagrangian:
 \begin{align}
			L_{\rho\beta}(\mathbf{x},\mathbf{z},\bm{\lambda})=&F(\mathbf{x})+H(\mathbf{z})-\langle(1-\rho\beta)\bm{\lambda},\mathbf{A}\mathbf{x}+\mathbf{B}\mathbf{z}-\mathbf{c}\rangle \notag \\
			&\qquad\qquad\qquad+\frac{\rho}{2} {\left\Vert \mathbf{A}\mathbf{x}+\mathbf{B}\mathbf{z}-\mathbf{c} \right\Vert }^2, \label{PerturbedAugmentedLagrangian}
		\end{align} 
        where  $\bm{\lambda}\in\mathbb{R}^p$ is the dual variable, $\rho>0$ is the penalty parameter, $\beta>0$ is the perturbation parameter, and satisfies $\rho\beta\in(0,1)$.
        We can observe that if $\beta=0$, \eqref{PerturbedAugmentedLagrangian} reduces to classical augmented Lagrangian $L_{\rho}(\mathbf{x},\mathbf{z},\bm{\lambda})=L_{\rho\beta}(\mathbf{x},\mathbf{z},\bm{\lambda})-\rho\beta\langle\bm{\lambda},\mathbf{A}\mathbf{x}+\mathbf{B}\mathbf{z}-\mathbf{c}\rangle$. 
    However, the necessity of positive $\beta$ is justified by the subsequent convergence analysis.

    To obtain closed-form solutions, we linearize the smooth components 
    \begin{equation*}
\begin{aligned}
F^0(\mathbf{x})\approx\langle\nabla  F^0(\mathbf{x}^k),\mathbf{x}-\mathbf{x}^k\rangle+F^0(\mathbf{x}^k),\\
H^0(\mathbf{z})\approx\langle\nabla  H^0(\mathbf{z}^k),\mathbf{z}-\mathbf{z}^k\rangle+H^0(\mathbf{z}^k),
\end{aligned}
\end{equation*}
  at the current iterates, and add proximal regularization terms 
    $$
	\frac{1}{2}{\Vert\mathbf{x}-\mathbf{x}^k\Vert}^2_{\mathbf{P}:=\tau_F\mathbf{I}-\rho\mathbf{A}^{\top}\mathbf{A}\succ\mathbf{0}}, \ \frac{1}{2}{\Vert\mathbf{z}-\mathbf{z}^k\Vert}^2_{\mathbf{Q}:=\tau_H\mathbf{I}-\rho\mathbf{B}^{\top}\mathbf{B}\succ\mathbf{0}},
	$$
    to the $\mathbf{x}$- and $\mathbf{z}$-subproblems, respectively.
These proximal terms are important to avoid the potential calculation of inverse matrices, which is particularly advantageous when $\mathbf{A}^{\top}\mathbf{A}$ and $\mathbf{B}^{\top}\mathbf{B}$ are rank-deficient or of high dimension \cite{wang2019distributed,bian2021stochastic}.
It yields the following explicit updates:
		\begin{align}
			\mathbf{x}^{k+1}\!=\!&\arg \min_{\mathbf{x}}\!\big\{\!\langle\nabla  F^0(\mathbf{x}^k),\mathbf{x}\rangle\!+\!F^1(\mathbf{x})\!+\!\frac{\rho}{2} {\left\Vert \mathbf{A}\mathbf{x}\!+\!\mathbf{B}\mathbf{z}^k\!-\!\mathbf{c} \right\Vert }^2 \notag\\
		&-\langle(1-\rho\beta)\bm{\lambda}^k,\mathbf{A}\mathbf{x}+\mathbf{B}\mathbf{z}^k-\mathbf{c}\rangle+\frac{1}{2}{\Vert\mathbf{x}-{\mathbf{x}}^k\Vert}^2_{\mathbf{P}} \big\}  \notag\\
        =&\textbf{prox}_{F^1}^{\tau_F}\{\tau_F^{-1}[-\nabla  F^0({\mathbf{x}}^k)+(\tau_F\mathbf{I}-\rho\mathbf{A}^{\top}\mathbf{A}){\mathbf{x}}^k \notag \\
& -\rho\mathbf{A}^{\top}\mathbf{B}\mathbf{z}^k+(1-\rho\beta)\mathbf{A}^{\top}\bm{\lambda}^k+\rho\mathbf{A}^{\top}\mathbf{c}] \}, \label{UpdateXProx}  \\ 
			\mathbf{z}^{k+1}\!=\!&\arg \min_{\mathbf{z}}  \big\{ \langle\nabla  H^0(\mathbf{z}^k),\!\mathbf{z}\rangle\!+\!H^1(\mathbf{z})\!+\!\frac{\rho}{2} {\left\Vert \mathbf{A}\mathbf{x}^{k+1}\!+\!\mathbf{B}\mathbf{z}-\!\mathbf{c} \right\Vert }^2 \notag\\
		&-\langle(1-\rho\beta)\bm{\lambda}^k, \mathbf{A}\mathbf{x}^{k+1}+\mathbf{B}\mathbf{z}-\mathbf{c}\rangle+\frac{1}{2}{\Vert\mathbf{z}-\mathbf{z}^k\Vert}^2_{\mathbf{Q}} \big\} \notag\\
=&\textbf{prox}_{H^1}^{\tau_H}\{\tau_H^{-1}[-\nabla  H^0(\mathbf{z}^k)+(\tau_H\mathbf{I}-\rho\mathbf{B}^{\top}\mathbf{B})\mathbf{z}^k \notag \\
& -\rho\mathbf{B}^{\top}\mathbf{A}\mathbf{x}^{k+1}+(1-\rho\beta)\mathbf{B}^{\top}\bm{\lambda}^k+\rho\mathbf{B}^{\top}\mathbf{c}] \}, \label{UpdateZProx}         \\
\bm{\lambda}^{k+1}=&(1-\rho\beta)\bm{\lambda}^{k}-\rho(\mathbf{A}\mathbf{x}^{k+1}+\mathbf{B}\mathbf{z}^{k+1}-\mathbf{c}), \label{ADMMla}
		\end{align} \label{ADMM}
where the proximal operator is defined as $\textbf{prox}_{g}^{\tau}(y)=\arg \min_x\{g(x)+\frac{\tau}{2}{\Vert x-y \Vert}^2\}$, with $\tau_F>\gamma_F$ and $\tau_H>\gamma_H$.
The detailed procedure are summarized in Algorithm 1.

\begin{algorithm}[t]\label{alg1}
\DontPrintSemicolon
\SetKwInOut{KwIn}{Initialization}
\SetKwData{Or}{\textbf{or}}
  \SetAlgoLined
  \KwIn {$\mathbf{x}^0$, $\mathbf{z}^0$ and $\bm{\lambda}^0$.}
  Set  $\tau_F\mathbf{I}\succ\rho\mathbf{A}^{\top}\mathbf{A}$, $\tau_H\mathbf{I}\succ\rho\mathbf{B}^{\top}\mathbf{B}$, $\tau_F>\gamma_F$, $\tau_H>\gamma_H$ \\
  \ and $\rho>0$, $\beta>0$, $0<\rho\beta<1$. \\
  \For{stopping criteria not satisfied}
  {
   Update the primal variable $\mathbf{x}^{k+1}$ by \eqref{UpdateXProx};\\
  Update the primal variable $\mathbf{z}^{k+1}$ by \eqref{UpdateZProx};\\
Update the dual variable $\bm{\lambda}^{k+1}$ by \eqref{ADMMla}.
  }
  \textbf{return} ${(\mathbf{x}^k,\mathbf{z}^k,\bm{\lambda}^k)}$.
  \caption{PPG-ADMM for \eqref{ClassicalADMMProblem}}
\end{algorithm}
	

    \begin{remark}
        {
      The perturbed dual update \eqref{ADMMla} can be interpreted as a dual ascent step derived from a Tikhonov-regularized Lagrangian with $\beta$-strong concavity in $\bm{\lambda}$ (i.e., $L_{\rho}(\mathbf{x},\mathbf{z},\bm{\lambda})-\frac{\beta}{2}{\Vert \bm{\lambda}\Vert}^2$, as in \cite{bernstein2019online,koshal2011multiuser}). 
      The negative quadratic term acts as ``damping", reducing  sensitivity in dual updates and enhancing numerical stability, particularly in ill-conditioned problems or under suboptimal parameter choices. 
      To provide a more intuitive illustration, we present a comparative example in Fig. \ref{fig:betaCase}, where the inclusion of perturbation term effectively prevents sustained oscillations, leading to faster convergence to a stable solution.}
       \begin{figure}[ht]
        \centering
        \includegraphics[width=0.9\linewidth]{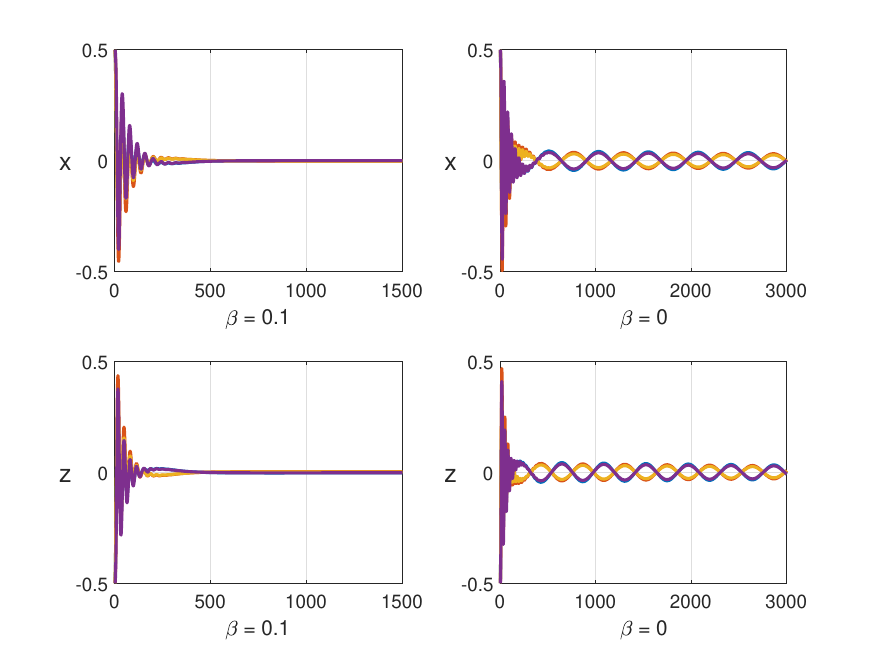}
        \vspace{-10pt}
        \caption{A simple example demonstrating the impact of $\beta$ on convergence behavior, where the figure depicts the iteration trajectories of all variable components.
        In this example, we use Algorithm 1 to solve  $\min_{\mathbf{x},\mathbf{z}}\Vert\mathbf{x} \Vert_1+\Vert\mathbf{z} \Vert_1$, $s.t.\ \mathbf{Ax}+\mathbf{Bz}=\mathbf{0}$, where $\mathbf{x}, \mathbf{z}\in\mathbb{R}^4$, and $\mathbf{A}$ and $\mathbf{B}$ are randomly generated rank-deficient matrices. 
        The optimal solution is $\mathbf{x}^*=\mathbf{z}^*=\mathbf{0}$.}\label{fig:betaCase}
        \vspace{-10pt}
    \end{figure}
    \end{remark}

\section{Convergence Analysis}\label{Section3}
In this section, we will analyze the convergence of the proposed PPG-ADMM framework. 
Although the specific problem formulations, constraints, and algorithms differ across existing relevant works \cite{li2015global,bian2021stochastic,boct2020proximal,guo2017convergence,yashtini2021multi,hajinezhad2018alternating,yashtini2022convergence,zhang2021extended,yang2017alternating,wang2018convergence,bot2019proximal,liu2019linearized,hong2016convergence,yang2022proximal,sun2019two,hajinezhad2019perturbed,wang2019global,wang2019distributed,melo2017iteration,zeng2024unified,huang2019faster,huang2018mini,zeng2024accelerated,chang2016asynchronous,hong2017distributed}, the convergence analysis typically follow a common  framework, which we also adopt.
The analysis begins with the construction of a Lyapunov function that is shown to decrease sufficiently and remain lower-bounded.
This step is often the most technically challenging, and it lays the foundation for the rest of the analysis.
Subsequently, the generated sequence is proved to be bounded and asymptotically regular.
Finally, it is established that the cluster point of the sequence is (near) the saddle point of the Lagrangian.

 One can utilize ADMM to address problem \eqref{ClassicalADMMProblem}:
\begin{subequations} \label{ClassicalADMM}
\begin{align}
\mathbf{x}^{k+1}=&\arg\min_{\mathbf{x}}L_{\rho}(\mathbf{x},\mathbf{z}^k,\bm{\lambda}^k), \label{ClassicalADMMxsubproblem} \\ 		
\mathbf{z}^{k+1}=&\arg\min_{\mathbf{z}}L_{\rho}(\mathbf{x}^{k+1},\mathbf{z},\bm{\lambda}^k), \label{ClassicalADMMzsubproblem}\\
\bm{\lambda}^{k+1}=&\bm{\lambda}^{k}-\rho(\mathbf{A}\mathbf{x}^{k+1}+\mathbf{B}\mathbf{z}^{k+1}-\mathbf{c}),\label{ClassicalADMMDualsubproblem}
\end{align} 
\end{subequations}
where $k\geqslant 0$ denotes the number of iterations, $\rho>0$ is the penalty parameter and $\bm{\lambda}\in\mathbb{R}^p$ is the dual variable.
The augmented Lagrangian in \eqref{ClassicalADMM} is defined as
\begin{align}
&L_{\rho}(\mathbf{x},\mathbf{z},\bm{\lambda}) \label{ClassicalADMMLagrangian}\\
=&F(\mathbf{x})+H(\mathbf{z})-\langle\bm{\lambda},\mathbf{A}\mathbf{x}+\mathbf{B}\mathbf{z}-\mathbf{c}\rangle+\frac{\rho}{2} {\left\Vert \mathbf{A}\mathbf{x}+\mathbf{B}\mathbf{z}-\mathbf{c} \right\Vert }^2. \notag
  \end{align}
  
Compared to the extensively studied convex optimization algorithms, relatively few algorithms are available for solving NCOPs.
Directly applying classic convex optimization methods to nonconvex problems is often infeasible, as certain common analytical techniques, such as variational inequalities, are generally not applicable in the nonconvex setting and thus fail to provide convergence guarantees \cite{wang2019global}.
In convergence analysis of some primal-dual algorithms for NCOPs, the crucial step is the construction of a Lyapunov function $\mathcal{P}$ \cite{wang2019global,yang2022proximal}, which must be both sufficiently decreasing and bounded from below with respect to the sequence $\{(\mathbf{x}^{k},\mathbf{z}^{k},\bm{\lambda}^{k})\}_{k\geqslant0}$ generated by algorithms, i.e.,
\begin{subequations}
\begin{align}
&\mathcal{P}(\mathbf{x}^{k+1},\mathbf{z}^{k+1},\bm{\lambda}^{k+1})-\mathcal{P}(\mathbf{x}^{k},\mathbf{z}^{k},\bm{\lambda}^{k}) \label{LyapunovCondition1}\\
\leqslant & -a_1{\Vert \mathbf{x}^{k+1}-\mathbf{x}^{k} \Vert}^2-a_2{\Vert \mathbf{z}^{k+1}-\mathbf{z}^{k} \Vert}^2-a_3{\Vert \bm{\lambda}^{k+1}-\bm{\lambda}^{k} \Vert}^2,  \notag\\
&\mathcal{P}(\mathbf{x}^{k},\mathbf{z}^{k},\bm{\lambda}^{k})>-\infty, \label{LyapunovCondition2}
\end{align}
\end{subequations}
where $a_1$, $a_2$ and $a_3$ are nonnegative coefficients.
Typically, such a Lyapunov function is constructed based on Lagrangian, as it facilitates establishing a connection between the cluster point of the generated sequence and the stationary solution of the optimization problem in subsequent analysis  \cite{yang2022proximal}.
However, directly employing the (augmented) Lagrangian as a Lyapunov function does not satisfy the requirements of \eqref{LyapunovCondition1}. 
Specifically, for 
\eqref{ClassicalADMMLagrangian}, we have
\begin{subequations}\label{successiveL}
\begin{align}
&L_{\rho}(\mathbf{x}^{k+1}\!,\!\mathbf{z}^{k},\bm{\lambda}^{k})\!-\!L_{\rho}(\mathbf{x}^{k},\mathbf{z}^{k},\bm{\lambda}^{k})\!\leqslant\! -a_1{\Vert \mathbf{x}^{k+1}\!-\!\mathbf{x}^{k} \Vert}^2\!,\! \label{successiveL1}\\
&L_{\rho}(\mathbf{x}^{k+1}\!,\mathbf{z}^{k+1},\bm{\lambda}^{k})\!-\!L_{\rho}(\mathbf{x}^{k+1}\!,\mathbf{z}^{k},\bm{\lambda}^{k})\!\leqslant\! -a_2{\Vert \mathbf{z}^{k+1}\!-\!\mathbf{z}^{k} \Vert}^2\!,\! \label{successiveL2}\\
    &\qquad L_{\rho}(\mathbf{x}^{k+1},\mathbf{z}^{k+1},\bm{\lambda}^{k+1})-L_{\rho}(\mathbf{x}^{k+1},\mathbf{z}^{k+1},\bm{\lambda}^{k}) \label{successiveL3}\\
    &\!=\!-\langle\bm{\lambda}^{k+1}-\bm{\lambda}^{k},\mathbf{A}\mathbf{x}^{k+1}+\mathbf{B}\mathbf{z}^{k+1}-\mathbf{c}\rangle\!\overset{\eqref{ClassicalADMMDualsubproblem}}{=}\!\frac{1}{\rho}{\Vert \bm{\lambda}^{k+1}-\bm{\lambda}^{k} \Vert}^2\!.\notag
\end{align}
\end{subequations}
The inequalities \eqref{successiveL1} and \eqref{successiveL2} follow from the optimality conditions of the primal variable updates, leveraging specific function properties such as smoothness or strong convexity.
However, the presence of positive coefficient $\rho$ in \eqref{successiveL3} contradicts the intention to construct a Lyapunov function satisfying \eqref{LyapunovCondition1}, after integrating \eqref{successiveL1}-\eqref{successiveL3}.
To address this, one can impose certain conditions to bound ${\Vert \bm{\lambda}^{k+1} - \bm{\lambda}^{k} \Vert}^2$ using successive differences of the primal variables, and then adjust the relevant coefficients accordingly to construct a valid Lyapunov function.
As a result, ADMM-based algorithms for solving NCOPs typically introduce the following two additional assumptions to ensure this control. 




In conventional ADMM convergence analysis, C1 and C2 are primarily introduced to facilitate the construction of a sufficiently decreasing Lyapunov function. 
We will first illustrate it in this section. 
To ensure that C1 and C2 hold, we assume $\mathrm{Im}(\mathbf{A}) \cup \{\mathbf{c}\} \subseteq \mathrm{Im}(\mathbf{B})$ and  ${H}^1 = 0$ in \eqref{ClassicalADMMProblem}.

Regarding the left-hand side (LHS) of \eqref{successiveL1}, we can directly deduce the following relation by using \eqref{ClassicalADMMxsubproblem}:
     \begin{equation}\label{FunctionLRelation1}
     L_{\rho}(\mathbf{x}^{k+1},\mathbf{z}^{k},\bm{\lambda}^{k})-L_{\rho}(\mathbf{x}^{k},\mathbf{z}^{k},\bm{\lambda}^{k})\leqslant 0.
 \end{equation}
 Then, for the LHS of \eqref{successiveL2}, we can further obtain
 \begin{align}
    &L_{\rho}(\mathbf{x}^{k+1},\mathbf{z}^{k+1},\bm{\lambda}^{k})-L_{\rho}(\mathbf{x}^{k+1},\mathbf{z}^{k},\bm{\lambda}^{k}) \label{FunctionLRelation2}\\
    =&H^0(\mathbf{z}^{k+1})-H^0(\mathbf{z}^{k})-\frac{\rho}{2} {\left\Vert \mathbf{z}^{k+1}-\mathbf{z}^{k}\right\Vert }^2_{\mathbf{B}^{\top}\mathbf{B}} \notag\\
    &+\langle\rho\mathbf{B}^{\top}(\mathbf{A}\mathbf{x}^{k+1}+\mathbf{B}\mathbf{z}^{k+1}-\mathbf{c})-\mathbf{B}^{\top}\bm{\lambda}^k,\mathbf{z}^{k+1}-\mathbf{z}^k\rangle \notag\\
    =&H^0(\mathbf{z}^{k+1})-H^0(\mathbf{z}^{k})-\langle\nabla  H^0(\mathbf{z}^{k}) ,\mathbf{z}^{k+1}-\mathbf{z}^k\rangle \notag\\
    &\!-\!\frac{\rho}{2} {\left\Vert \mathbf{z}^{k+1}\!-\!\mathbf{z}^{k}\right\Vert }^2_{\mathbf{B}^{\top}\mathbf{B}}\!-\!\langle\nabla H^0(\mathbf{z}^{k+1})\!-\!\nabla H^0(\mathbf{z}^{k}) ,\mathbf{z}^{k+1}\!-\!\mathbf{z}^k\rangle \notag\\
        \leqslant&-\frac{1}{2}{\Vert\mathbf{z}^{k+1}-\mathbf{z}^k\Vert}^2_{\rho\mathbf{B}^{\top}\mathbf{B}-3L_H\mathbf{I}}, \notag
\end{align}
where the last equation is derived from the optimality condition of \eqref{ClassicalADMMzsubproblem}, i.e.,
\begin{equation}
\begin{aligned}\label{OptimalityConditionSmooth}
\nabla H^0(\mathbf{z}^{k+1})&=\mathbf{B}^{\top}\bm{\lambda}^k-\rho\mathbf{B}^{\top}(\mathbf{A}\mathbf{x}^{k+1}+\mathbf{B}\mathbf{z}^{k+1}-\mathbf{c}) \\
&\overset{\eqref{ClassicalADMMDualsubproblem}}{=}\mathbf{B}^{\top}\bm{\lambda}^{k+1},
\end{aligned}
\end{equation}
while the last inequality relies on the smoothness of $H^0$.
Subsequently, we demonstrate that both C1 and C2 are related to bounding the term on the right-hand side (RHS) of \eqref{successiveL3}.
In addition to ensuring the feasibility, C1 plays a role in the initial step of bounding ${\Vert \bm{\lambda}^{k+1}-\bm{\lambda}^{k} \Vert}^2$.
Specifically, through the use of the dual variable update \eqref{ClassicalADMMDualsubproblem}, it can be inferred that $(\bm{\lambda}^{k+1}-\bm{\lambda}^{k})\in\mathrm{Im}(\mathbf{B})$, thereby implying that \cite[Lemma 3]{wang2019global}
\begin{equation}\label{dualVariable}
    {\Vert \bm{\lambda}^{k+1}-\bm{\lambda}^{k} \Vert}^2\leqslant \frac{1}{\sigma_{\mathbf{B}\mathbf{B}^{\top}}^+}{\Vert \mathbf{B}^{\top}(\bm{\lambda}^{k+1}-\bm{\lambda}^{k}) \Vert}^2,
\end{equation}
where $\sigma_{\mathbf{B}\mathbf{B}^{\top}}^+$ denotes the smallest positive eigenvalue of ${\mathbf{B}\mathbf{B}^{\top}}$.
Note that \eqref{dualVariable} also holds for rank-deficient matrices, and thus $\mathbf{B}$ is not required to be full rank or the identity matrix.
Next, we can use the successive difference of the primal variable $\mathbf{z}$ to provide an upper bound for the RHS of \eqref{dualVariable}:
\begin{align}
&{\Vert \mathbf{B}^{\top}(\bm{\lambda}^{k+1}-\bm{\lambda}^{k}) \Vert}\overset{\eqref{ClassicalADMMDualsubproblem}}{=}{\Vert \rho\mathbf{B}^{\top}(\mathbf{A}\mathbf{x}^{k+1}+\mathbf{B}\mathbf{z}^{k+1}-\mathbf{c}) \Vert} \notag\\
\overset{\eqref{OptimalityConditionSmooth}}{=}&{\Vert \nabla H^0(\mathbf{z}^{k+1})-\nabla H^0(\mathbf{z}^{k}) \Vert}\leqslant L_H{\Vert \mathbf{z}^{k+1}-\mathbf{z}^{k} \Vert}.\label{OptimalityConditionIn}
\end{align}
The property of smooth function is used in this inequality. 
Similar results can be found in \cite{hajinezhad2018alternating,zhang2021extended,guo2017convergence,wang2019distributed,hong2016convergence,chang2016asynchronous,hong2017distributed}.
Additionally, as described in \cite{sun2019two}, as $k\rightarrow\infty$, the residual ${\Vert \mathbf{A}\mathbf{x}^{k+1}+\mathbf{B}\mathbf{z}^{k+1}-\mathbf{c} \Vert}$ may not be eliminated. 
However, \eqref{OptimalityConditionIn} indicates that it can be controlled by ${\Vert \mathbf{z}^{k+1}-\mathbf{z}^{k} \Vert}$, further emphasizing the necessity of C2.
Combining \eqref{FunctionLRelation1}, \eqref{FunctionLRelation2}, \eqref{dualVariable} and \eqref{OptimalityConditionIn}, we have
\begin{align*}
    &L_{\rho}(\mathbf{x}^{k+1},\mathbf{z}^{k+1},\bm{\lambda}^{k+1})-L_{\rho}(\mathbf{x}^{k},\mathbf{z}^{k},\bm{\lambda}^{k})  \\
    \leqslant&-\frac{1}{2}{\Vert\mathbf{z}^{k+1}-\mathbf{z}^k\Vert}^2_{\rho\mathbf{B}^{\top}\mathbf{B}}+\big(\frac{L_H^2}{\rho\sigma_{\mathbf{B}\mathbf{B}^{\top}}^+}+\frac{3L_H}{2}\big){\Vert \mathbf{z}^{k+1}-\mathbf{z}^k \Vert}^2.
\end{align*}
By adjusting $\rho$, a sufficiently decreasing function is obtained.

Apart from the above analysis, the condition $\mathrm{Im}(\mathbf{A}) \cup \{\mathbf{c}\} \subseteq \mathrm{Im}(\mathbf{B})$ ensures the feasibility of the nonconvex optimization problem and guarantees that, when the Gauss-Seidel type algorithm converges to some $\mathbf{x}^*$ and $\mathbf{z}^*$, the constraint $\mathbf{A}\mathbf{x}^*+\mathbf{B}\mathbf{z}^*=\mathbf{c}$ holds \cite{wang2019global,sun2019two,melo2017iteration}.
Moreover, this range space inclusion assumption is commonly used to establish the boundedness of dual variables and the convergence of the primal residuals. 
This naturally raises the question of whether an alternative condition can be established to bound the successive differences of dual variables and to construct a decreasing Lyapunov function, thereby eliminating the need for C2.
Under this relaxed setting, if we can further prove the boundedness of all variables, the column space inclusion assumption in C1 can also be dropped, and feasibility alone would suffice (i.e., Assumption \ref{AssumptionMatrix}).
Fortunately, our subsequent analysis reveals that the perturbation introduced in PPG-ADMM makes this relaxation possible.


We introduce the following function $\mathcal{T}$, which incorporates the perturbed augmented Lagrangian \eqref{PerturbedAugmentedLagrangian}, two proximal terms from the primal variable updates, and an additional term related to the dual variable:
 \begin{equation*}
	\begin{aligned}
		&\mathcal{T}(\mathbf{x},\mathbf{z},\bm{\lambda},{\mathbf{x}}{'},\mathbf{z}{'})=L_{\rho\beta}(\mathbf{x},\mathbf{z},\bm{\lambda})+\frac{1}{2}{\Vert\mathbf{x}-{\mathbf{x}}'\Vert}^2_{\mathbf{P}}+\frac{1}{2}{\Vert\mathbf{z}-\mathbf{z}'\Vert}^2_{\mathbf{Q}}  \\
  &\qquad\qquad\qquad\qquad-\frac{\beta}{2}(1-\rho\beta){\Vert\bm{\lambda}\Vert}^2.
	\end{aligned}
 \end{equation*}
	\begin{proposition} \label{PropositionStronglyConvex}
		Function $\mathcal{T}(\mathbf{x},\mathbf{z},\bm{\lambda},{\mathbf{x}}{'},\mathbf{z}{'})$ is $(\tau_F-\gamma_F-L_F)$-strongly convex with respect to $\mathbf{x}$ and $(\tau_H-\gamma_H-L_H)$-strongly convex with respect to $\mathbf{z}$.
	\end{proposition}
	\begin{IEEEproof}
We begin by analyzing the smooth component of $\mathcal{T}$ with respect to $\mathbf{x}$.
  For any $\overline{\mathbf{x}}, \underline{\mathbf{x}}\in\mathbb{R}^n$, we have
  \begin{equation*}
		\begin{aligned}
			&\langle\nabla _{\mathbf{x}} [\mathcal{T}(\overline{\mathbf{x}},\mathbf{z},\bm{\lambda},{\mathbf{x}}{'},\mathbf{z}{'})-F^1(\overline{\mathbf{x}})]	 \\
			&-\nabla _{\mathbf{x}} [\mathcal{T}(\underline{\mathbf{x}},\mathbf{z},\bm{\lambda},{\mathbf{x}}{'},\mathbf{z}{'})-F^1(\underline{\mathbf{x}})], \overline{\mathbf{x}}-\underline{\mathbf{x}}\rangle \\
			=&\langle\nabla\!F^0(\overline{\mathbf{x}})\!-\!\nabla\!  F^0(\underline{\mathbf{x}})\!+\!\tau_F(\overline{\mathbf{x}}\!-\!\underline{\mathbf{x}}), \overline{\mathbf{x}}-\underline{\mathbf{x}}\rangle\!\geqslant\!(\tau_F\!-\!L_F){\Vert\overline{\mathbf{x}}-\underline{\mathbf{x}}\Vert}^2\!,
		\end{aligned}
  \end{equation*}
		where the inequality follows from the $L_F$-smoothness of $F^0$.
  It can be inferred that $\mathcal{T}(\mathbf{x},\mathbf{z},\bm{\lambda},{\mathbf{x}}{'},\mathbf{z}{'})-F^1({\mathbf{x}})$ is $(\tau_F-L_F)$-strongly convex with respect to $\mathbf{x}$.
		Since $F^1(\mathbf{x})$ is $\gamma_F$-weakly convex, it follows that $\mathcal{T}$ is $(\tau_F-\gamma_F-L_F)$-strongly convex with respect to $\mathbf{x}$. 
  Similarly, 
        we can deduce that the function $\mathcal{T}$ is $(\tau_H-\gamma_H-L_H)$-strongly convex with respect to $\mathbf{z}$.
	\end{IEEEproof}


The optimality conditions corresponding to the updates \eqref{UpdateXProx} and \eqref{UpdateZProx}  are given by:
\begin{equation}\label{OptimalityConditionX}
\begin{aligned}
&\nabla  F^0({\mathbf{x}}^k)+\mathbf{v}^{F^1}_{\mathbf{x}^{k+1}}+\rho\mathbf{A}^{\top}  
(\mathbf{A}\mathbf{x}^{k+1}+\mathbf{B}\mathbf{z}^{k}-\mathbf{c}) \\
&\quad-(1-\rho\beta)\mathbf{A}^{\top}\bm{\lambda}^k+\mathbf{P}(\mathbf{x}^{k+1}-{\mathbf{x}}^{k})=\mathbf{0}, 
    \end{aligned}
    \end{equation}
\begin{equation}\label{OptimalityConditionZ}
    \begin{aligned}
&\nabla  H^0(\mathbf{z}^k)+\mathbf{v}_{\mathbf{z}^{k+1}}^{H^1}+\rho\mathbf{B}^{\top}(\mathbf{A}\mathbf{x}^{k+1}+\mathbf{B}\mathbf{z}^{k+1}-\mathbf{c}) \\
&\quad-(1-\rho\beta)\mathbf{B}^{\top}\bm{\lambda}^k+\mathbf{Q}(\mathbf{z}^{k+1}-\mathbf{z}^{k})=\mathbf{0}, 
    \end{aligned}
    \end{equation}
    where $\mathbf{v}^{F^1}_{\mathbf{x}^{k+1}}\in\partial F^1(\mathbf{x}^{k+1})$ and $\mathbf{v}_{\mathbf{z}^{k+1}}^{H^1}\in\partial H^1(\mathbf{z}^{k+1})$ are subgradients of the respective nonsmooth components.
By leveraging these optimality conditions together with the updates in Algorithm 1, we can derive key relationships about the successive differences of the primal and dual variables.
	\begin{proposition}\label{PropositionOptimalityCondition}
Suppose that Assumptions \ref{AssumptionMatrix}, \ref{AssumptionFunction} hold, and let $\{(\mathbf{x}^{k},\mathbf{z}^{k},\bm{\lambda}^{k})\}_{k\geqslant0}$ be the sequence generated by Algorithm 1.
   Then, for $k\geqslant 1$, the following inequality holds:
		\begin{align}\label{Proposition2}
		&{\Vert\mathbf{x}^{k+1}-\mathbf{x}^{k}\Vert}^2_{L_F\mathbf{I}+\mathbf{P}}-{\Vert\mathbf{x}^{k}-\mathbf{x}^{k-1}\Vert}^2_{L_F\mathbf{I}+\mathbf{P}} \notag\\
		&+{\Vert\mathbf{z}^{k+1}\!-\mathbf{z}^{k}\Vert}^2_{L_H\mathbf{I}+\mathbf{Q}+2\rho\mathbf{B}^{\top}\mathbf{B}}\!-\!{\Vert\mathbf{z}^{k}-\mathbf{z}^{k-1}\Vert}^2_{L_H\mathbf{I}+\mathbf{Q}+2\rho\mathbf{B}^{\top}\mathbf{B}} \notag\\
		&+\frac{1-\rho\beta}{\rho}{\Vert\bm{\lambda}^{k+1}-\bm{\lambda}^{k}\Vert}^2-\frac{1-\rho\beta}{\rho}{\Vert\bm{\lambda}^{k}-\bm{\lambda}^{k-1}\Vert}^2 \\
		\leqslant&{\Vert\mathbf{x}^{k+1}-\mathbf{x}^k\Vert}^2_{2(L_F+\gamma_F)\mathbf{I}+\rho\mathbf{A}^{\top}\mathbf{A}} \notag\\
		&+{\Vert\mathbf{z}^{k+1}-\mathbf{z}^k\Vert}^2_{2(L_H+\gamma_H)\mathbf{I}+4\rho\mathbf{B}^{\top}\mathbf{B}}-2\beta{\Vert\bm{\lambda}^{k+1}-\bm{\lambda}^{k}\Vert}^2.\notag
		\end{align}
	\end{proposition}
	\begin{IEEEproof}
 First, the inner product of \eqref{OptimalityConditionX} and $\mathbf{x}-\mathbf{x}^{k+1}$ is
		\begin{align} \label{OptimalityConditionXk+1}
  &\langle\nabla  F^0({\mathbf{x}}^k)+\rho\mathbf{A}^{\top}  
(\mathbf{A}\mathbf{x}^{k+1}+\mathbf{B}\mathbf{z}^{k}-\mathbf{c}) \notag\\
		&-(1-\rho\beta)\mathbf{A}^{\top}\bm{\lambda}^k+\mathbf{P}(\mathbf{x}^{k+1}-{\mathbf{x}}^{k}), \mathbf{x}-\mathbf{x}^{k+1} \rangle \notag\\
		\overset{\eqref{ADMMla}}{=}&\langle\nabla  F^0(\mathbf{x}^k)-\rho\mathbf{A}^{\top}\mathbf{B}(\mathbf{z}^{k+1}-\mathbf{z}^{k})+\mathbf{P}(\mathbf{x}^{k+1}-\mathbf{x}^{k}) \notag\\
		&-\mathbf{A}^{\top}\bm{\lambda}^{k+1}, \mathbf{x}-\mathbf{x}^{k+1}\rangle=-\langle\mathbf{v}^{F^1}_{\mathbf{x}^{k+1}}, \mathbf{x}-\mathbf{x}^{k+1}\rangle.
		\end{align}
		Analogously, the optimality condition for $\mathbf{x}^k$ gives:
\begin{align} \label{OptimalityConditionXk}
  &\langle\nabla  F^0({\mathbf{x}}^{k-1})-\rho\mathbf{A}^{\top}\mathbf{B}(\mathbf{z}^{k}-\mathbf{z}^{k-1})+\mathbf{P}(\mathbf{x}^{k}-\mathbf{x}^{k-1}) \notag\\
&-\mathbf{A}^{\top}\bm{\lambda}^{k}, \mathbf{x}-\mathbf{x}^{k} \rangle =-\langle\mathbf{v}^{F^1}_{\mathbf{x}^{k}}, \mathbf{x}-\mathbf{x}^{k}\rangle.
		\end{align}
  By substituting $\mathbf{x}=\mathbf{x}^k$ into \eqref{OptimalityConditionXk+1} and $\mathbf{x}=\mathbf{x}^{k+1}$ into \eqref{OptimalityConditionXk}, and subsequently adding them together, we obtain
		\begin{align}\label{OptimalityConditionXResult}
			&\langle\mathbf{x}^{k}-\mathbf{x}^{k+1}, \nabla  F^0(\mathbf{x}^k)-\nabla  F^0(\mathbf{x}^{k-1})\rangle \notag\\
			&+\langle\mathbf{x}^{k}-\mathbf{x}^{k+1}, -\rho\mathbf{A}^{\top}\mathbf{B}[(\mathbf{z}^{k+1}-\mathbf{z}^{k})-(\mathbf{z}^{k}-\mathbf{z}^{k-1})]\rangle \notag\\
			&+\langle\mathbf{x}^{k}-\mathbf{x}^{k+1}, \mathbf{P}[(\mathbf{x}^{k+1}-\mathbf{x}^{k})-(\mathbf{x}^{k}-\mathbf{x}^{k-1})]\rangle \notag\\
			&+\langle\mathbf{x}^{k}-\mathbf{x}^{k+1}, -\mathbf{A}^{\top}(\bm{\lambda}^{k+1}-\bm{\lambda}^{k})\rangle \notag\\
			=&\langle\mathbf{v}^{F^1}_{\mathbf{x}^{k+1}}-\mathbf{v}^{F^1}_{\mathbf{x}^{k}}, \mathbf{x}^{k+1}-\mathbf{x}^{k}\rangle{\geqslant}-\gamma_F{\Vert\mathbf{x}^{k+1}-\mathbf{x}^{k}\Vert}^2, 
		\end{align}
		where the inequality holds due to the weakly convexity of $F^1$.
		Next, we will analyze each term on the left-hand side (LHS) of \eqref{OptimalityConditionXResult} in detail.
		 By using Young's inequality and the smoothness assumption, the first term satisfies
		\begin{align}\label{OptimalityConditionProof1}
		&\langle\mathbf{x}^{k}-\mathbf{x}^{k+1}, \nabla  F^0(\mathbf{x}^k)-\nabla  F^0(\mathbf{x}^{k-1})\rangle \notag\\
		\leqslant&\frac{L_F}{2}{\Vert \mathbf{x}^{k+1}-\mathbf{x}^{k}\Vert}^2+\frac{1}{2L_F}{\Vert \nabla  F^0(\mathbf{x}^k)-\nabla  F^0(\mathbf{x}^{k-1})\Vert}^2 \notag\\
		\leqslant&\frac{L_F}{2}{\Vert \mathbf{x}^{k+1}-\mathbf{x}^{k}\Vert}^2+\frac{L_F}{2}{\Vert \mathbf{x}^{k}-\mathbf{x}^{k-1}\Vert}^2.
		\end{align}
		Similarly, by applying Young's inequality to the second term of \eqref{OptimalityConditionXResult}, we can deduce that
		\begin{align}\label{OptimalityConditionProof2}
		&\rho\langle\mathbf{A}(\mathbf{x}^{k+1}-\mathbf{x}^{k}), \mathbf{B}[(\mathbf{z}^{k+1}-\mathbf{z}^{k})-(\mathbf{z}^{k}-\mathbf{z}^{k-1})] \\
		\leqslant&\frac{1}{2}{\Vert \mathbf{x}^{k+1}\!-\!\mathbf{x}^{k}\Vert}^2_{\rho\mathbf{A}^{\top}\mathbf{A}}\!+\!\frac{1}{2}{\Vert (\mathbf{z}^{k+1}\!-\!\mathbf{z}^{k})\!-\!(\mathbf{z}^{k}\!-\!\mathbf{z}^{k-1})\Vert}^2_{\rho\mathbf{B}^{\top}\mathbf{B}} \notag \\
		\leqslant&\frac{1}{2}{\Vert \mathbf{x}^{k+1}\!-\!\mathbf{x}^{k}\Vert}^2_{\rho\mathbf{A}^{\top}\mathbf{A}}\!+\!{\Vert \mathbf{z}^{k+1}\!-\!\mathbf{z}^{k}\Vert}^2_{\rho\mathbf{B}^{\top}\mathbf{B}}\!+\!{\Vert \mathbf{z}^{k}\!-\!\mathbf{z}^{k-1}\Vert}^2_{\rho\mathbf{B}^{\top}\mathbf{B}}. \notag
		\end{align}
		By using $\langle -x, x-y\rangle=-\frac{1}{2}({\Vert x-y \Vert}^2+{\Vert x \Vert}^2-{\Vert y \Vert}^2)$, the third term of \eqref{OptimalityConditionXResult} can be transformed into
		\begin{align}\label{OptimalityConditionProof3}
			&\langle\mathbf{x}^{k}-\mathbf{x}^{k+1}, \mathbf{P}[(\mathbf{x}^{k+1}-\mathbf{x}^{k})-(\mathbf{x}^{k}-\mathbf{x}^{k-1})]\rangle \notag\\
			=&-\frac{1}{2}{\Vert(\mathbf{x}^{k+1}-\mathbf{x}^{k})-(\mathbf{x}^{k}-\mathbf{x}^{k-1})\Vert}^2_{\mathbf{P}}-\frac{1}{2}{\Vert\mathbf{x}^{k+1}-\mathbf{x}^{k}\Vert}^2_{\mathbf{P}} \notag\\
			&+\frac{1}{2}{\Vert\mathbf{x}^{k}-\mathbf{x}^{k-1}\Vert}^2_{\mathbf{P}}.
		\end{align}

        Now, we can derive the conditions for $\mathbf{z}^{k+1}$ and $\mathbf{z}^k$, following the same logic as \eqref{OptimalityConditionXk+1} and \eqref{OptimalityConditionXk}:
		\begin{align}
		&\langle \nabla  H^0(\mathbf{z}^k)+\rho\mathbf{B}^{\top}(\mathbf{A}\mathbf{x}^{k+1}+\mathbf{B}\mathbf{z}^{k+1}-\mathbf{c}) \notag\\
		&-(1-\rho\beta)\mathbf{B}^{\top}\bm{\lambda}^k+\mathbf{Q}(\mathbf{z}^{k+1}-\mathbf{z}^{k}), \mathbf{z}-\mathbf{z}^{k+1}\rangle \notag\\
		\overset{\eqref{ADMMla}}{=}&\langle \nabla  H^0(\mathbf{z}^k)-\mathbf{B}^{\top}\bm{\lambda}^{k+1}+\mathbf{Q}(\mathbf{z}^{k+1}-\mathbf{z}^{k}), \mathbf{z}-\mathbf{z}^{k+1}\rangle \notag\\
		=&-\langle\mathbf{v}_{\mathbf{z}^{k+1}}^{H^1}, \mathbf{z}-\mathbf{z}^{k+1}\rangle, \label{OptimalityConditionZk+1}
		\end{align} 
  		\begin{align}
  &\langle \nabla  H^0(\mathbf{z}^{k-1})-\mathbf{B}^{\top}\bm{\lambda}^{k}+\mathbf{Q}(\mathbf{z}^{k}-\mathbf{z}^{k-1}), \mathbf{z}-\mathbf{z}^{k}\rangle \notag\\
		=&-\langle \mathbf{v}_{\mathbf{z}^{k}}^{H^1}, \mathbf{z}-\mathbf{z}^{k}\rangle. \label{OptimalityConditionZk}
		\end{align}
		We set $\mathbf{z}=\mathbf{z}^k$ and $\mathbf{z}=\mathbf{z}^{k+1}$ in \eqref{OptimalityConditionZk+1} and \eqref{OptimalityConditionZk}, respectively, and add them together, which yields a result similar to \eqref{OptimalityConditionXResult}:
	\begin{align}\label{OptimalityConditionZResult}
	&\langle\mathbf{z}^k-\mathbf{z}^{k+1}, \nabla  H^0(\mathbf{z}^k)-\nabla  H^0(\mathbf{z}^{k-1})\rangle \notag\\
	&+\langle\mathbf{z}^k-\mathbf{z}^{k+1}, \mathbf{Q}[(\mathbf{z}^{k+1}-\mathbf{z}^{k})-(\mathbf{z}^{k}-\mathbf{z}^{k-1})]\rangle \notag\\
	&+\langle\mathbf{z}^k-\mathbf{z}^{k+1},-\mathbf{B}^{\top}(\bm{\lambda}^{k+1}-\bm{\lambda}^{k})\rangle \notag\\
=&\langle\mathbf{v}^{H^1}_{\mathbf{z}^{k+1}}-\mathbf{v}^{H^1}_{\mathbf{z}^{k}}, \mathbf{z}^{k+1}-\mathbf{z}^{k}\rangle{\geqslant} -\gamma_H{\Vert\mathbf{z}^{k+1}-\mathbf{z}^{k}\Vert}^2.
		\end{align}
  Similar to \eqref{OptimalityConditionProof1} and \eqref{OptimalityConditionProof3}, the first and second terms on the LHS of \eqref{OptimalityConditionZResult} can be bounded as		\begin{align}\label{OptimalityConditionProof4}
			&\langle\mathbf{z}^k-\mathbf{z}^{k+1}, \nabla  H^0(\mathbf{z}^k)-\nabla  H^0(\mathbf{z}^{k-1})\rangle \notag\\
			\leqslant& \frac{L_H}{2}{\Vert\mathbf{z}^{k+1}-\mathbf{z}^k\Vert}^2+\frac{L_H}{2}{\Vert\mathbf{z}^{k}-\mathbf{z}^{k-1}\Vert}^2,
		\end{align}
\begin{align}\label{OptimalityConditionProof5}
&\langle\mathbf{z}^k-\mathbf{z}^{k+1}, \mathbf{Q}[(\mathbf{z}^{k+1}-\mathbf{z}^{k})-(\mathbf{z}^{k}-\mathbf{z}^{k-1})]\rangle \notag\\
=&-\frac{1}{2}{\Vert(\mathbf{z}^{k+1}-\mathbf{z}^{k})-(\mathbf{z}^{k}-\mathbf{z}^{k-1})\Vert}^2_{\mathbf{Q}}-\frac{1}{2}{\Vert\mathbf{z}^{k+1}-\mathbf{z}^{k}\Vert}^2_{\mathbf{Q}} \notag\\
			&+\frac{1}{2}{\Vert\mathbf{z}^{k}-\mathbf{z}^{k-1}\Vert}^2_{\mathbf{Q}}.
		\end{align}
		Next, by combining the dual-related terms on the LHS of \eqref{OptimalityConditionXResult} and \eqref{OptimalityConditionZResult}, we have 
\begin{align}
&\langle\bm{\lambda}^{k+1}-\bm{\lambda}^{k}, (\mathbf{A}\mathbf{x}^{k+1}+\mathbf{B}\mathbf{z}^{k+1})-(\mathbf{A}\mathbf{x}^{k}+\mathbf{B}\mathbf{z}^{k})\rangle \notag\\
\overset{\eqref{ADMMla}}{=}&-\frac{1}{\rho}{\Vert\bm{\lambda}^{k+1}-\bm{\lambda}^{k}\Vert}^2+\frac{1-\rho\beta}{\rho}\langle\bm{\lambda}^{k+1}-\bm{\lambda}^{k}, \bm{\lambda}^{k}-\bm{\lambda}^{k-1}\rangle \notag\\
=&\frac{1-\rho\beta}{2\rho}{\Vert\bm{\lambda}^{k}-\bm{\lambda}^{k-1}\Vert}^2-\frac{1+\rho\beta}{2\rho}{\Vert\bm{\lambda}^{k+1}-\bm{\lambda}^{k}\Vert}^2 \notag\\
&-\frac{1-\rho\beta}{2\rho}{\Vert(\bm{\lambda}^{k+1}-\bm{\lambda}^{k})-(\bm{\lambda}^{k}-\bm{\lambda}^{k-1})\Vert}^2. \label{OptimalityConditionProof6}
		\end{align}
		Finally, by combining \eqref{OptimalityConditionXResult}-\eqref{OptimalityConditionProof3}, \eqref{OptimalityConditionZResult}-\eqref{OptimalityConditionProof6} and omitting the non-positive terms on the right-hand side (RHS), \eqref{Proposition2} can be obtained.
	\end{IEEEproof}

 Next, we analyze the successive difference of the function $\mathcal{T}$ with respect to $\{\bm{\omega}^{k+1}:=(\mathbf{x}^{k+1},\mathbf{z}^{k+1},\bm{\lambda}^{k+1},\mathbf{x}^k,\mathbf{z}^k)\}_{k\geqslant0}$ generated by PPG-ADMM.

\begin{proposition}\label{PropositionRelationT}
    Suppose that Assumptions \ref{AssumptionMatrix} and \ref{AssumptionFunction} hold.
    For $k\geqslant 1$, $\mathcal{T}(\bm{\omega}^{k+1})$ and $\mathcal{T}(\bm{\omega}^{k})$ satisfy the following inequality:
		\begin{align}
			&\mathcal{T}(\bm{\omega}^{k+1})-\mathcal{T}(\bm{\omega}^{k})\leqslant -\frac{\tau_F-\gamma_F-3L_F}{2}{\Vert\mathbf{x}^{k+1}-\mathbf{x}^k\Vert}^2\notag\\
			&\!-\!\frac{1}{2}{\Vert\mathbf{x}^k\!-\!\mathbf{x}^{k-1}\Vert}^2_{\mathbf{P}}\!-\!\frac{\tau_H-\gamma_H-3L_H}{2}{\Vert\mathbf{z}^{k+1}\!-\!\mathbf{z}^k\Vert}^2 \notag\\ 
   &\!-\!\frac{1}{2}{\Vert\mathbf{z}^k\!-\!\mathbf{z}^{k-1}\Vert}^2_{\mathbf{Q}}\!+\!\frac{(1-\rho\beta)(2-\rho\beta)}{2\rho}{\Vert\bm{\lambda}^{k+1}\!-\!\bm{\lambda}^{k}\Vert}^2.\label{FunctionTRelation}
		\end{align}
	\end{proposition}
	\begin{IEEEproof}
 By utilizing the strong convexity of $\mathcal{T}$ with respect to $\mathbf{x}$ in Proposition \ref{PropositionStronglyConvex}, it yields the following inequality:
\begin{align}\label{FunctionTProof1}
	&\mathcal{T}(\mathbf{x}^{k+1},\mathbf{z}^{k},\bm{\lambda}^{k},{\mathbf{x}}^k,\mathbf{z}^k)\!-\!\mathcal{T}({\mathbf{x}}^{k},\mathbf{z}^{k},\bm{\lambda}^{k},{\mathbf{x}}^k,\mathbf{z}^k) \\
	{\leqslant}& \langle \mathbf{v}^{\mathcal{T}}_{\mathbf{x}^{k+1}}, \mathbf{x}^{k+1}-{\mathbf{x}}^{k}\rangle-\frac{\tau_F-\gamma_F-L_F}{2}{\Vert\mathbf{x}^{k+1}-{\mathbf{x}}^k\Vert}^2 \notag\\
	=&\langle \nabla F^0(\mathbf{x}^{k+1})+\mathbf{v}^{F^1}_{\mathbf{x}^{k+1}}-(1-\rho\beta)\mathbf{A}^{\top}\bm{\lambda}^k+\mathbf{P}(\mathbf{x}^{k+1}-{\mathbf{x}}^{k}) \notag\\
 &+\rho\mathbf{A}^{\top}(\mathbf{A}\mathbf{x}^{k+1}+\mathbf{B}\mathbf{z}^{k}-\mathbf{c}), \mathbf{x}^{k+1}-{\mathbf{x}}^k \rangle \notag\\
	&-\frac{\tau_F-\gamma_F-L_F}{2}{\Vert\mathbf{x}^{k+1}-{\mathbf{x}}^k\Vert}^2 \notag\\
	\overset{\eqref{OptimalityConditionX}}{=}& \langle \nabla F^0(\mathbf{x}^{k+1})-\nabla F^0({\mathbf{x}}^{k}), \mathbf{x}^{k+1}-{\mathbf{x}}^k \rangle \notag\\
     &-\frac{\tau_F-\gamma_F-L_F}{2}{\Vert\mathbf{x}^{k+1}-{\mathbf{x}}^k\Vert}^2 \notag \\
     \leqslant& -\frac{\tau_F-\gamma_F-3L_F}{2}{\Vert\mathbf{x}^{k+1}-{\mathbf{x}}^k\Vert}^2, \notag
		\end{align}
		where  $\mathbf{v}^{\mathcal{T}}_{\mathbf{x}^{k+1}}\in\partial_\mathbf{x} \mathcal{T}(\mathbf{x}^{k+1},\mathbf{z}^{k},\bm{\lambda}^{k},\mathbf{x}^k,\mathbf{z}^k)$.
		By applying a similar argument to $\mathbf{z}$, we obtain:
\begin{align}\label{FunctionTProof2}
&\mathcal{T}(\mathbf{x}^{k+1},\mathbf{z}^{k+1},\bm{\lambda}^{k},{\mathbf{x}}^k,\mathbf{z}^k)-\!\mathcal{T}(\mathbf{x}^{k+1},\mathbf{z}^{k},\bm{\lambda}^{k},{\mathbf{x}}^k,\mathbf{z}^k) \\
\leqslant&-\frac{\tau_H-\gamma_H-3L_H}{2}{\Vert\mathbf{z}^{k+1}-\mathbf{z}^k\Vert}^2, \notag
\end{align}
Next, we focus on the dual variable $\bm{\lambda}$:
\begin{align}\label{FunctionTProof3}
&\mathcal{T}(\bm{\omega}^{k+1})-\mathcal{T}({\mathbf{x}}^{k+1},\mathbf{z}^{k+1},\bm{\lambda}^{k},{\mathbf{x}}^k,\mathbf{z}^k) \\
=& -(1\!-\!\rho\beta)[\langle\bm{\lambda}^{k+1}\!-\!\bm{\lambda}^{k}, \mathbf{A}{\mathbf{x}}^{k+1}+\mathbf{B}\mathbf{z}^{k+1}-\mathbf{c}\rangle+\frac{\beta}{2}{\Vert\bm{\lambda}^{k+1}\Vert}^2 \notag\\
&-\frac{\beta}{2}{\Vert\bm{\lambda}^{k}\Vert}^2]=\frac{(1-\rho\beta)(2-\rho\beta)}{2\rho}{\Vert\bm{\lambda}^{k+1}-\bm{\lambda}^{k}\Vert}^2,\notag
\end{align}
where the last equality is derived from \eqref{ADMMla} and $\langle x, x-y\rangle=\frac{1}{2}({\Vert x-y \Vert}^2+{\Vert x \Vert}^2-{\Vert y \Vert}^2)$.	
Finally, for the last two parameters of $\mathcal{T}$, we can directly obtain
\begin{align}\label{FunctionTProof4}
&\mathcal{T}(\mathbf{x}^{k},\mathbf{z}^{k},\bm{\lambda}^{k},\mathbf{x}^k,\mathbf{z}^k)-\mathcal{T}(\bm{\omega}^{k}) \notag\\
=&-\frac{1}{2}{\Vert\mathbf{x}^k-\mathbf{x}^{k-1}\Vert}^2_{\mathbf{P}}-\frac{1}{2}{\Vert\mathbf{z}^k-\mathbf{z}^{k-1}\Vert}^2_{\mathbf{Q}}. 
\end{align}
Combining \eqref{FunctionTProof1}, \eqref{FunctionTProof2}, \eqref{FunctionTProof3} and \eqref{FunctionTProof4} establishes  \eqref{FunctionTRelation}.
	\end{IEEEproof}

Next, we construct a new function $\mathcal{P}$, which is composed of $\mathcal{T}$ and three terms on the LHS of \eqref{Proposition2}:
\begin{align*}
&\mathcal{P}(\bm{\omega}^{k+1})=\mathcal{T}(\bm{\omega}^{k+1})+d\left[{\Vert\mathbf{x}^{k+1}-\mathbf{x}^{k}\Vert}^2_{L_F\mathbf{I}+\mathbf{P}} \right.\\
  &\left.+{\Vert\mathbf{z}^{k+1}-\mathbf{z}^{k}\Vert}^2_{L_H\mathbf{I}+\mathbf{Q}+2\rho\mathbf{B}^{\top}\mathbf{B}}+\frac{1-\rho\beta}{\rho}{\Vert\bm{\lambda}^{k+1}-\bm{\lambda}^{k}\Vert}^2 \right],
	\end{align*}
  where $d>0$ is an adjustable coefficient.
  Since both matrices $\mathbf{P}$ and $\mathbf{Q}$ are positive definite, and $\rho,{L_F},L_H>0$, it follows that $L_F\mathbf{I}+\mathbf{P}$ and $L_H\mathbf{I}+\mathbf{Q}+2\rho\mathbf{B}^{\top}\mathbf{B}$ are also positive definite matrices.
  We make the following assumption about the parameters in the function $\mathcal{P}$.

\begin{assumption}\label{AssumptionCoefficients}
     The parameters in the Lyapunov function $\mathcal{P}$ need to satisfy
     \begin{subequations}
         \begin{align}
        &\tau_F\mathbf{I}\succ2d\rho\mathbf{A}^{\top}\mathbf{A}+[(4d+3)L_F+(4d+1)\gamma_F]\mathbf{I}, \label{AssumptionCoefficients1}\\
        &\tau_H\mathbf{I}\succ8d\rho\mathbf{B}^{\top}\mathbf{B}+[(4d+3)L_H+(4d+1)\gamma_H]\mathbf{I}, \label{AssumptionCoefficients2}\\   
        &d>\frac{(1-\rho\beta)(2-\rho\beta)}{4\rho\beta}>0, \label{AssumptionCoefficients3}\\
        &\rho>0,\ \beta>0,\ 0<\rho\beta<1,\label{AssumptionCoefficients4}\\ &\tau_F\mathbf{I}\succ\rho\mathbf{A}^{\top}\mathbf{A},  \ \tau_H\mathbf{I}\succ\rho\mathbf{B}^{\top}\mathbf{B}, \ \tau_F>\gamma_F, \ \tau_H>\gamma_H. \label{AssumptionCoefficients5}
         \end{align}
     \end{subequations}
 \end{assumption}
Note that \eqref{AssumptionCoefficients4} and \eqref{AssumptionCoefficients5} have already been assumed in the development of this algorithm.
Since $\tau_F$, $\tau_H$ and $d$ are adjustable parameters, \eqref{AssumptionCoefficients1}-\eqref{AssumptionCoefficients3} can be easily satisfied.

\begin{proposition}\label{Theorem1}
Suppose that Assumptions \ref{AssumptionMatrix}, \ref{AssumptionFunction} and \ref{AssumptionCoefficients} hold.
Then, the sequence $\{(\mathbf{x}^{k},\mathbf{z}^{k},\bm{\lambda}^{k})\}_{k\geqslant0}$ generated by Algorithm 1 is bounded, and the Lyapunov function $\mathcal{P}$ is sufficiently decreasing and bounded from below. 
The successive differences of the iterates vanish asymptotically, i.e., $\Vert \mathbf{x}^{k+1}-\mathbf{x}^{k} \Vert\rightarrow 0$, $\Vert \mathbf{z}^{k+1}-\mathbf{z}^{k} \Vert\rightarrow 0$, ${\Vert \bm{\lambda}^{k+1}-\bm{\lambda}^{k} \Vert}\rightarrow 0$, as $k\rightarrow +\infty$.
\end{proposition}
\begin{IEEEproof}
   (i) \textbf{Sufficient Descent of $\mathcal{P}$}. 
   By combining inequalities \eqref{Proposition2} and \eqref{FunctionTRelation}, we obtain the following relationship between $\mathcal{P}(\bm{\omega}^{k+1})$ and $\mathcal{P}(\bm{\omega}^{k})$:
	\begin{align}\label{FunctionPRelation}
		&\mathcal{P}(\bm{\omega}^{k+1})-\mathcal{P}(\bm{\omega}^{k}) \notag\\
		\leqslant&-\frac{1}{2}{\Vert\mathbf{x}^{k+1}-\mathbf{x}^k\Vert}^2_{[\tau_F-(4d+1)\gamma_F-(4d+3)L_F]\mathbf{I}-2d\rho\mathbf{A}^{\top}\mathbf{A}} \notag\\
		&-\frac{1}{2}{\Vert\mathbf{z}^{k+1}-\mathbf{z}^k\Vert}^2_{[\tau_H-(4d+1)\gamma_H-(4d+3)L_H]\mathbf{I}-8d\rho\mathbf{B}^{\top}\mathbf{B}} \notag\\
  &-\frac{1}{2}{\Vert\mathbf{x}^k-\mathbf{x}^{k-1}\Vert}^2_{\mathbf{P}}-\frac{1}{2}{\Vert\mathbf{z}^k-\mathbf{z}^{k-1}\Vert}^2_{\mathbf{Q}} \notag\\
		&-\left[2d\beta-\frac{(1-\rho\beta)(2-\rho\beta)}{2\rho}\right]{\Vert\bm{\lambda}^{k+1}-\bm{\lambda}^{k}\Vert}^2. 
	\end{align}
 Note that the conditions \eqref{AssumptionCoefficients1}, \eqref{AssumptionCoefficients2} and \eqref{AssumptionCoefficients3} in Assumption \ref{AssumptionCoefficients} are imposed to ensure that all terms on the RHS of \eqref{FunctionPRelation} are negative.
 This directly implies that $\mathcal{P}$ is sufficiently decreasing with respect to $\{\bm{\omega}^{k+1}\}_{k\geqslant0}$.
 A direct consequence is that $\mathcal{P}(\bm{\omega}^{k+1})$ is uniformly upper-bounded by its initial value, i.e., $\mathcal{P}(\bm{\omega}^{k+1}) \leqslant \mathcal{P}(\bm{\omega}^{1}), \forall k\geqslant0$.

(ii) \textbf{Boundedness of the Iterate Sequence}.
From the definition of $\mathcal{P}$, since $F$ and $H$ are bounded from below, and the remaining terms are non-negative, there exists a constant $\hat{\mathcal{P}}$ associated with $\mathcal{P}(\bm{\omega}^{1})$ such that
\begin{align}\label{FunctionPUpper}
    \hat{\mathcal{P}}\geqslant&-\langle(1-\rho\beta)\bm{\lambda}^{k+1},\mathbf{A}\mathbf{x}^{k+1}+\mathbf{B}\mathbf{z}^{k+1}-\mathbf{c}+\frac{\beta}{2}\bm{\lambda}^{k+1}\rangle  \notag\\
    \overset{\eqref{ADMMla}}{=}&\frac{(1-\rho\beta)^2}{\rho}\langle\bm{\lambda}^{k+1}, \bm{\lambda}^{k+1}-\bm{\lambda}^{k}\rangle+\frac{\beta(1-\rho\beta)}{2}{\Vert \bm{\lambda}^{k+1} \Vert}^2 \notag \\
    =&\frac{(1-\rho\beta)^2}{2\rho}\left[{\Vert \bm{\lambda}^{k+1}-\bm{\lambda}^{k} \Vert}^2+{\Vert \bm{\lambda}^{k+1} \Vert}^2-{\Vert \bm{\lambda}^{k} \Vert}^2\right] \notag \\
    &+\frac{\beta(1-\rho\beta)}{2}{\Vert \bm{\lambda}^{k+1} \Vert}^2. 
\end{align}
We now proceed by induction.
Assume that the
constant $\hat{\mathcal{P}}$ is chosen such that
$$
\frac{\beta(1-\rho\beta)}{2}{\Vert \bm{\lambda}^{k} \Vert}^2\leqslant\hat{\mathcal{P}},\ {\forall}k\leqslant K.
$$
If ${\Vert \bm{\lambda}^{K+1} \Vert}\geqslant{\Vert \bm{\lambda}^{K} \Vert}$, then by \eqref{FunctionPUpper}, it has 
\begin{align*}
    \hat{\mathcal{P}}&\geqslant\frac{(1-\rho\beta)^2}{2\rho}\left[{\Vert \bm{\lambda}^{K+1} \Vert}^2-{\Vert \bm{\lambda}^{K} \Vert}^2\right]+\frac{\beta(1-\rho\beta)}{2}{\Vert \bm{\lambda}^{K+1} \Vert}^2\\
    &\geqslant \frac{\beta(1-\rho\beta)}{2}{\Vert \bm{\lambda}^{K+1} \Vert}^2.
\end{align*}
If instead ${\Vert \bm{\lambda}^{K+1} \Vert}<{\Vert \bm{\lambda}^{K} \Vert}$, then it directly follows that
$$
\frac{\beta(1-\rho\beta)}{2}{\Vert \bm{\lambda}^{K+1} \Vert}^2<\frac{\beta(1-\rho\beta)}{2}{\Vert \bm{\lambda}^{K} \Vert}^2\leqslant\hat{\mathcal{P}}.
$$
Therefore, in either case, $\{\bm{\lambda}^k\}_{k \geqslant 0}$ remains bounded.
Next, by combining \eqref{FunctionPUpper} and the definition of $\mathcal{P}(\bm{\omega}^{k+1})$, we can conclude that
\begin{align}
    &F(\mathbf{x}^{k+1})+H(\mathbf{z}^{k+1})+\frac{1}{2}{\Vert\mathbf{x}^{k+1}-\mathbf{x}^{k}\Vert}^2_{2dL_F\mathbf{I}+(2d+1)\mathbf{P}} \notag\\
    &+\frac{1}{2}{\Vert\mathbf{z}^{k+1}-\mathbf{z}^{k}\Vert}^2_{2dL_H\mathbf{I}+4d\rho\mathbf{B}^{\top}\mathbf{B}+(2d+1)\mathbf{Q}} \notag\\
    &+\frac{\rho}{2} {\left\Vert \mathbf{A}\mathbf{x}^{k+1}+\mathbf{B}\mathbf{z}^{k+1}-\mathbf{c} \right\Vert }^2 \notag\\
    &+\frac{(2d+2-\rho\beta)(1-\rho\beta)}{2\rho}{\Vert\bm{\lambda}^{k+1}-\bm{\lambda}^{k}\Vert}^2<+\infty, \label{FunctionPBound2}
\end{align}
where the terms related to ${\Vert \bm{\lambda} \Vert}$ is omitted since they are already known to be bounded.
By Assumption \ref{AssumptionCoefficients}, all involved matrices are positive definite and all coefficients are strictly positive.
Consequently, we conclude that $\{\Vert\mathbf{x}^{k+1}-\mathbf{x}^{k}\Vert \}_{k\geqslant0}$, $\{\Vert\mathbf{z}^{k+1}-\mathbf{z}^{k}\Vert \}_{k\geqslant0}$ and $\{\Vert\bm{\lambda}^{k+1}-\bm{\lambda}^{k}\Vert\}_{k\geqslant0}$ are bounded.
Moreover, since both $F$ and $H$ are coercive by Assumption \ref{AssumptionFunction} \ref{AssumptionFunctioni}, and $F(\mathbf{x}^{k+1})$, $H(\mathbf{z}^{k+1})$ are upper-bounded, this implies that the sequences $\{\mathbf{x}^k\}_{k\geqslant0}$ and $\{\mathbf{z}^k\}_{k\geqslant0}$ must also be bounded.

(iii) \textbf{Asymptotic Regularity and Lower Bound of $\mathcal{P}$}.
We first establish that $\mathcal{P}$ is bounded from below. 
This is a direct consequence of the boundedness of the iterate sequence. 
Specifically, the inner product term $-\langle(1-\rho\beta)\bm{\lambda}^{k+1},\mathbf{A}\mathbf{x}^{k+1}+\mathbf{B}\mathbf{z}^{k+1}-\mathbf{c}+\frac{\beta}{2}\bm{\lambda}^{k+1}\rangle$ within $\mathcal{P}(\bm{\omega}^{k+1})$ can be shown to be bounded below by leveraging the boundedness of $\{\bm{\lambda}^k\}_{k \geqslant 0}$ and \eqref{FunctionPUpper}, which relates it to bounded quadratic terms of $\bm{\lambda}$. 
As all other terms in the definition of $\mathcal{P}$ are either non-negative by construction or bounded below by Assumption \ref{AssumptionFunction} \ref{AssumptionFunctioni}, it follows that the entire sequence $\{\mathcal{P}(\bm{\omega}^{k+1})\}_{k\geqslant0}$ is bounded from below.
Now, since the sequence is both monotonically decreasing and bounded from below, it must converge to a finite limit. 
Recalling \eqref{FunctionPBound2}, all non-negative terms must converge to $0$.
This proves the asymptotic regularity claims, and the proof is thus completed.
\end{IEEEproof}

The necessity of the condition $\beta \neq 0$ is further underscored by \eqref{FunctionPRelation}.
Specifically, when $\beta = 0$, the coefficient of the last term on the RHS of \eqref{FunctionPRelation} becomes a positive constant $1/\rho$, which cannot be rendered negative by adjusting $d$.
Under such scenario, where C1 and C2 cannot simultaneously hold, constructing a valid Lyapunov function becomes nearly infeasible without imposing another assumptions.

Subsequently, we provide the definition of (Approximate) Karush–Kuhn–Tucker (AKKT/KKT) point for problem \eqref{ClassicalADMMProblem}.
Both KKT point and AKKT point involve \eqref{AKKT1} and \eqref{AKKT2}; the distinction lies in the slight violation of the equality constraint.
\begin{definition}\label{def:AKKT}
Consider the following inequalities with $\epsilon>0$
    \begin{subequations}\label{AKKT}
      \begin{align}
&\mathrm{dist}(\mathbf{0},\nabla F^0(\mathbf{x}^{*})+\partial F^1(\mathbf{x}^{*})-\mathbf{A}^{\top}\bm{\lambda}^*)\leqslant\epsilon, \label{AKKT1}\\
&\mathrm{dist}(\mathbf{0},\nabla H^0(\mathbf{z}^{*})+\partial H^1(\mathbf{z}^{*})-\mathbf{B}^{\top}\bm{\lambda}^*)\leqslant\epsilon, \label{AKKT2}\\
&\|\mathbf{A}\mathbf{x}^{*}+\mathbf{B}\mathbf{z}^{*}-\mathbf{c}+\beta\bm{\lambda}^*\|\leqslant\epsilon,\label{AKKT3} \\
&\|\mathbf{A}\mathbf{x}^{*}+\mathbf{B}\mathbf{z}^{*}-\mathbf{c}\|\leqslant\epsilon. \label{AKKT4}
      \end{align}
  \end{subequations}
 If a point $(\mathbf{x}^{*},\mathbf{z}^{*},\bm{\lambda}^{*})$ satisfies \eqref{AKKT1}, \eqref{AKKT2} and \eqref{AKKT3}, it is referred to as an $\epsilon$-AKKT point (or approximate stationary point).
 It is called an $\epsilon$-KKT point if it satisfies \eqref{AKKT1}, \eqref{AKKT2} and \eqref{AKKT4}.
 If $\epsilon=0$, the point $(\mathbf{x}^{*},\mathbf{z}^{*},\bm{\lambda}^{*})$ is said to satisfy the AKKT/KKT conditions \cite{zhang2021online} for \eqref{ClassicalADMMProblem}.
\end{definition}

\begin{theorem}\label{Theorem2}
    Suppose that Assumptions  \ref{AssumptionMatrix}, \ref{AssumptionFunction} and \ref{AssumptionCoefficients} hold. 
Then, for $1\leqslant j\leqslant K$, there exists a point ${(\mathbf{x}^{j+1},\mathbf{z}^{j+1},\bm{\lambda}^{j+1})}$, such that if $K={\zeta^2}/{\epsilon^2}$ for some positive constant $\zeta$, it will be an $\epsilon$-AKKT point.
Moreover, if $\beta=1/K$ and $K={(\zeta+\sqrt{2\hat{\mathcal{P}}})^2}/{\epsilon^2}+\rho$, with $\hat{\mathcal{P}}$ being a finite constant used to bound $\{\bm{\lambda}^{k}\}_{k\geqslant0}$, it will be an $\epsilon$-KKT point.   
\end{theorem}
  \begin{IEEEproof}
  By Proposition \ref{Theorem1},  the sequence $\{(\mathbf{x}^{k}, \mathbf{z}^{k}, \bm{\lambda}^{k})\}_{k\geqslant0}$ generated by Algorithm 1 is bounded and will converge to a cluster point.
From the optimality conditions \eqref{OptimalityConditionX}, \eqref{OptimalityConditionZ}, and the dual update \eqref{ADMMla}, we obtain
\begin{align}
 &\nabla F^0(\mathbf{x}^{k+1})\!-\!\nabla  F^0(\mathbf{x}^k)\!-\!\mathbf{P}(\mathbf{x}^{k+1}\!-\!\mathbf{x}^{k})\!+\!\rho\mathbf{A}^{\top}\mathbf{B}(\mathbf{z}^{k+1}\!-\!\mathbf{z}^{k}) \notag\\
 &\qquad\qquad\quad\in\nabla  F^0(\mathbf{x}^{k+1})+\partial F^1(\mathbf{x}^{k+1})-\mathbf{A}^{\top}\bm{\lambda}^{k+1},\label{AKKTC1}\\
 &\nabla  H^0(\mathbf{z}^{k+1})-\nabla  H^0(\mathbf{z}^k)-\mathbf{Q}(\mathbf{z}^{k+1}-\mathbf{z}^{k}) \notag\\
 &\qquad\qquad\quad\in \nabla  H^0(\mathbf{z}^{k+1})+\partial H^1(\mathbf{z}^{k+1})-\mathbf{B}^{\top}\bm{\lambda}^{k+1}.\label{AKKTC2}
\end{align}
Summing \eqref{FunctionPRelation} over $k=1,\dots, K$, there exists a finite upper bound $\bar{\mathcal{P}}=2\mathcal{P}(\bm{\omega}^{1})-2\mathcal{P}(\bm{\omega}^{K+1})+{\Vert\mathbf{x}^{K+1}-\mathbf{x}^K\Vert}^2_{\mathbf{P}}-{\Vert\mathbf{x}^{1}-\mathbf{x}^0\Vert}^2_{\mathbf{P}}+{\Vert\mathbf{z}^{K+1}-\mathbf{z}^K\Vert}^2_{\mathbf{Q}}-{\Vert\mathbf{z}^{1}-\mathbf{z}^0\Vert}^2_{\mathbf{Q}}$ that satisfies
\begin{equation*}
   	\begin{aligned}
\bar{\mathcal{P}}\geqslant&\sum_{k=1}^{K}{\Vert\mathbf{x}^{k+1}-\mathbf{x}^k\Vert}^2_{[\tau_F-(4d+1)\gamma_F-(4d+3)L_F]\mathbf{I}-2d\rho\mathbf{A}^{\top}\mathbf{A}+\mathbf{P}} \notag\\
&+\sum_{k=1}^{K}{\Vert\mathbf{z}^{k+1}-\mathbf{z}^k\Vert}^2_{[\tau_H-(4d+1)\gamma_H-(4d+3)L_H]\mathbf{I}-8d\rho\mathbf{B}^{\top}\mathbf{B}+\mathbf{Q}} \notag\\
&+[4d\beta-{(1-\rho\beta)(2-\rho\beta)}/{\rho}]\sum_{k=1}^K{\Vert\bm{\lambda}^{k+1}-\bm{\lambda}^{k}\Vert}^2\geqslant0.
	\end{aligned}
\end{equation*}
It means that there exists $1\leqslant j\leqslant K$ such that
 $$
    {\Vert\mathbf{x}^{j+1}-\mathbf{x}^j\Vert}\leqslant\sqrt{\frac{\bar{\mathcal{P}}}{\sigma_1^{\min}K}}, \  {\Vert\mathbf{z}^{j+1}-\mathbf{z}^j\Vert}\leqslant\sqrt{\frac{\bar{\mathcal{P}}}{\sigma_2^{\min}K}},
$$
$$
    {\Vert\bm{\lambda}^{j+1}-\bm{\lambda}^j\Vert}\leqslant\sqrt{\frac{\rho\bar{\mathcal{P}}}{[4d\rho\beta-(1-\rho\beta)(2-\rho\beta)]K}},
$$
where $\sigma_1^{\min}$ and $\sigma_2^{\min}$ are the smallest eigenvalues of the positive definite matrices $[\tau_F-(4d+1)\gamma_F-(4d+3)L_F]\mathbf{I}-2d\rho\mathbf{A}^{\top}\mathbf{A}+\mathbf{P}$ and $[\tau_H-(4d+1)\gamma_H-(4d+3)L_H]\mathbf{I}-8d\rho\mathbf{B}^{\top}\mathbf{B}+\mathbf{Q}$, respectively.
Therefore, for the LHS of \eqref{AKKTC1} and \eqref{AKKTC2}, it has 
{\begin{align}
&\|\nabla\! F^0(\mathbf{x}^{j+1})\!-\!\nabla\!  F^0(\mathbf{x}^j)\!-\!\mathbf{P}(\mathbf{x}^{j+1}\!-\!\mathbf{x}^{j})\!+\!\rho\mathbf{A}^{\top}\!\mathbf{B}(\mathbf{z}^{j+1}\!-\!\mathbf{z}^{j})\| \notag\\
\leqslant&\|\nabla F^0(\mathbf{x}^{j+1})-\nabla  F^0(\mathbf{x}^j)\|+\|\mathbf{P}(\mathbf{x}^{j+1}-\mathbf{x}^{j})\| \notag\\
&+\rho\|\mathbf{A}^{\top}\mathbf{B}(\mathbf{z}^{j+1}-\mathbf{z}^{j})\| \notag\\
\leqslant& (L_F+\|\mathbf{P}\|)\|\mathbf{x}^{j+1}-\mathbf{x}^{j}\|+\rho\|\mathbf{A}^{\top}\mathbf{B}\|\|\mathbf{z}^{j+1}-\mathbf{z}^{j}\| \notag\\
\leqslant& {\frac{\zeta_1:=(L_F+\|\mathbf{P}\|)\sqrt{\bar{\mathcal{P}}/\sigma_1^{\min}}\!+\!\rho\|\mathbf{A}^{\top}\mathbf{B}\|\sqrt{\bar{\mathcal{P}}/\sigma_2^{\min}}}{\sqrt{K}}}, \label{t2e1}
\end{align}
\begin{align}
&\|\nabla  H^0(\mathbf{z}^{j+1})-\nabla  H^0(\mathbf{z}^j)-\mathbf{Q}(\mathbf{z}^{j+1}-\mathbf{z}^{j})\| \notag\\
    \leqslant&\|\nabla  H^0(\mathbf{z}^{j+1})-\nabla  H^0(\mathbf{z}^j)\|+\|\mathbf{Q}\|\|\mathbf{z}^{j+1}-\mathbf{z}^{j}\| \notag\\
    \leqslant&\frac{\zeta_2:=(L_H+\|\mathbf{Q}\|)\sqrt{\bar{\mathcal{P}}/\sigma_2^{\min}}}{\sqrt{K}}.\label{t2e2}
\end{align}}
Moreover, by using \eqref{ADMMla}, it also has
\begin{align}
&\|\mathbf{A}\mathbf{x}^{j+1}+\mathbf{B}\mathbf{z}^{j+1}-\mathbf{c}+\beta\bm{\lambda}^{j+1}\|=\frac{1+\rho\beta}{\rho}\|\bm{\lambda}^{j+1}-\bm{\lambda}^{j}\| \notag\\
\leqslant&{\frac{\zeta_3:=(1+\rho\beta)\sqrt{\bar{\mathcal{P}}/[4d\rho^2\beta-\rho(1-\rho\beta)(2-\rho\beta)]}}{\sqrt{K}}}.\label{t2e3}
\end{align}
Defining ${\zeta\!:=\!\max\{\zeta_1, \zeta_2, \zeta_3\}}/{\sqrt{K}}\!=\!\epsilon$, we deduce $K\!=\!\zeta/\epsilon^2$.
Then, by \eqref{t2e1}, \eqref{t2e2} and \eqref{t2e3}, the point ${(\mathbf{x}^{j+1},\mathbf{z}^{j+1},\bm{\lambda}^{j+1})}$ satisfies the conditions of an $\epsilon$-AKKT point.

 To further obtain an $\epsilon$-KKT point, we consider the bound on the equality constraint:
\begin{align}
&\|\mathbf{A}\mathbf{x}^{j+1}+\mathbf{B}\mathbf{z}^{j+1}-\mathbf{c}\| \label{t2e4}\\
\leqslant&\|\mathbf{A}\mathbf{x}^{j+1}+\mathbf{B}\mathbf{z}^{j+1}-\mathbf{c}+\beta\bm{\lambda}^{j+1}\|+\beta\|\bm{\lambda}^{j+1}\| \notag\\
=&\frac{1+\rho\beta}{\rho}\|\bm{\lambda}^{j+1}-\bm{\lambda}^{j}\|+\beta\|\bm{\lambda}^{j+1}\| \leqslant\frac{\zeta}{\sqrt{K}}+\sqrt{\frac{2\beta\hat{\mathcal{P}}}{1-\rho\beta}}.\notag
\end{align}
We note that the feasibility of problem \eqref{ClassicalADMMProblem} has been assumed, ensuring the existence of points satisfying the equality constraint.
If we set $\beta = 1/K$, then 
$$\frac{\zeta}{\sqrt{K}}+\sqrt{\frac{2\beta\hat{\mathcal{P}}}{1-\rho\beta}}=\frac{\zeta}{\sqrt{K}}+\sqrt{\frac{2\hat{\mathcal{P}}}{K-\rho}}\leqslant\frac{\zeta+\sqrt{2\hat{\mathcal{P}}}}{\sqrt{K-\rho}}=\epsilon.$$
Therefore, we can set $K={(\zeta+\sqrt{2\hat{\mathcal{P}}})^2}/{\epsilon^2}+\rho$.
Combining \eqref{t2e1}, \eqref{t2e2} and \eqref{t2e4}, the point ${(\mathbf{x}^{j+1},\mathbf{z}^{j+1},\bm{\lambda}^{j+1})}$ satisfies the definition of an $\epsilon$-KKT point.
  \end{IEEEproof}
  Theorem \ref{Theorem2} demonstrates that PPG-ADMM exhibits a sublinear convergence rate of $\mathcal{O}(1/\sqrt{K})$ and an iteration complexity of $\mathcal{O}(1/\epsilon^2)$, which aligns with the known lower complexity bounds for nonconvex optimization problems \cite{chen2021communication}.

\begin{figure}[htbp]
\centering
\vspace{-15pt}
\includegraphics[width=0.8\linewidth]{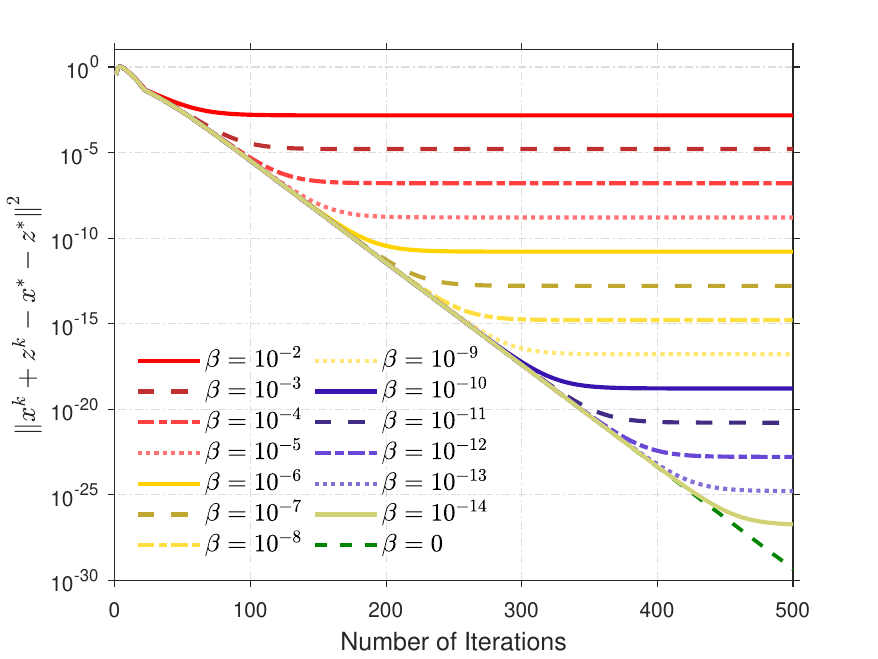}
\vspace{-5pt}
\caption{A simple case to evaluate the accuracy of PPG-ADMM. 
The objective function is $x^3+2(x-1)^2+\phi(x)+z^3+2(z-1)^2+\mathcal{M}_{1,1}(z)$, subject to $x+z=0$, where $\phi$ is the indicator function of the set $\mathcal{C}=\{x\in\mathbb{R}|-2\leqslant x \leqslant 2\}$, $\mathcal{M}_{1,1}(z)$ is $1$-weakly convex MCP.
The optimal solution is $x^*=z^*=0$.}\label{fig:remark_test}
\vspace{-10pt}
\end{figure}

\begin{remark}
Theoretically, the choice of $\beta$ directly affects the constraint violation of PPG-ADMM. 
To illustrate this, we present a simple test example in Fig. \ref{fig:remark_test} to evaluate the accuracy. 
As observed, larger values of $\beta$ lead to faster convergence to a stationary point but incur greater constraint violations, whereas smaller values achieve the opposite. 
This highlights an inherent trade-off.
In particular, even when $\beta = 0$, the algorithm still converges in this example, albeit in the absence of any theoretical guarantees.
\end{remark}

\section{Application: Decentralized Partial Consensus}

Distributed optimization relies on the collaboration of agents in the network to achieve the global objective.
 To describe the network structure, we utilize the undirected graph $\mathcal{G}=\{\mathcal{V},\mathcal{E}\}$, where $\mathcal{V}$ is the set of agents, and $\mathcal{E}\subseteq\mathcal{V}\times\mathcal{V}$ represents the set of undirected edges connecting agents within $\mathcal{V}$. 
 For instance, $(i,j)\in\mathcal{E}$ indicates a direct connection between agents $i$ and $j$.
The neighbor set of $i\in\mathcal{V}$ is represented by $\mathcal{N}_i=\{j\in\mathcal{V}|(i,j)\in\mathcal{E},i\neq j\}$.  
 Each agent $i$ aims to minimize a local objective function $f_i(x_i)$ while satisfying coupled constraints with its neighbors. 
 We formulate this as the partial consensus problem  on $\mathcal{G}$ \cite{koppel2017proximity,hajinezhad2019perturbed}:
\begin{equation}\label{EdgeProblem}
\begin{aligned}
    	\min  &\quad  \sum_{i=1}^{N}{f_i}({x}_i)=\sum_{i=1}^{N}f_i^0({x}_i)+\sum_{i=1}^{N}f_i^1({x}_i), \\
     s.t. &\quad A_{ij}x_i+A_{ji}x_j\in \mathcal{C}_{(i,j)}, \ (i,j)\in\mathcal{E},
\end{aligned}
\end{equation}
	where $x_i\in\mathbb{R}^{n_i}$, $x_j\in\mathbb{R}^{n_j}$,  $A_{ij}\in\mathbb{R}^{l_{(i,j)}\times n_i}$ and $A_{ji}\in\mathbb{R}^{l_{(i,j)}\times n_j}$ are linear mapping matrices, $\mathcal{C}_{(i,j)}$ is a convex set defining the feasible relationship between agents. 
    This formulation generalizes the standard consensus problem (where $x_i=x_j$) and accommodates more complex coordination tasks.

To apply the proposed PPG-ADMM, we introduce an auxiliary variable $\mathbf{z}$ that concatenates the variables $z_{(i,j)}$ for all edges, allowing us to rewrite the constraint in \eqref{EdgeProblem} as an equality constraint $A_{ij}x_i+A_{ji}x_j+z_{(i,j)}=0$ with $-z_{(i,j)}\in \mathcal{C}_{(i,j)}$.
Then, the problem is reformulated into the standard composite form
\begin{equation*}
			\min \  F(\mathbf{x})+H(\mathbf{z}), \quad s.t. \ {\mathbf{N}}\mathbf{x}+\mathbf{z}=\mathbf{0}, 
	\end{equation*}
	where $\mathbf{x}=[x_1;\cdots;x_N]\in\mathbb{R}^{\sum_{i=1}^Nn_i}$, $\mathbf{z}\in\mathbb{R}^{\sum_{(i,j)\in\mathcal{E}}l_{(i,j)}}$, $F(\mathbf{x})=F^0(\mathbf{x})+F^1(\mathbf{x})=\sum_{i=1}^Nf_i(x_i)$.
    The operator $\mathbf{N}\in\mathbb{R}^{\sum_{(i,j)\in\mathcal{E}}l_{(i,j)}\times\sum_{i=1}^Nn_i}$ is block-stacked from $N_{(i,j)}:{x}\mapsto A_{ij}x_i+A_{ji}x_j$.
    Moreover,  $H(\mathbf{z})=\sum_{(i,j)\in\mathcal{E}}h_{(i,j)}(z_{(i,j)})$ is the sum of indicator functions for the sets $\mathcal{C}_{(i,j)}$.
Thanks to the separability of $F$ and $H$, and the block-sparse structure of $\mathbf{N}$, iterations \eqref{UpdateXProx}-\eqref{ADMMla} naturally decompose into local operations.
Specifically, the local updates for agent $i\in\mathcal{V}$ are as follows:
\begin{subequations}\label{DistributedEdgeProblemUpdate}
\begin{align}
x_i^{k+1}=&\textbf{prox}_{f_i^1}^{\tau_F}\{\tau^{-1}_1[-\nabla  f_i^0({x}^k_i)+(\tau_F-\rho\sum_{j\in\mathcal{N}_i}A_{ij}^{\top}A_{ij}){x}^k_i  \notag\\
&-\rho\sum_{j\in\mathcal{N}_i}A_{ij}^{\top}A_{ji}{x}^k_j-\rho\sum_{j\in\mathcal{N}_i}A_{ij}^{\top}z_{(i,j),i}^{k} \notag\\
&+(1-\rho\beta)\sum_{j\in\mathcal{N}_i}A_{ij}^{\top}\lambda_{(i,j),i}^{k}] \},\label{DistributedEdgeProblemUpdate1}\\
z_{(i,j),i}^{k+1}=&\delta_{\mathcal{C}_{(i,j)}}\{\tau_H^{-1}[(\tau_H-\rho)z_{(i,j),i}^{k}-\rho(A_{ij}x_i^{k+1}+A_{ji}x_j^{k+1}) \notag\\
&+(1-\rho\beta)\lambda_{(i,j),i}^{k}]\}, \  \forall j\in\mathcal{N}_i, \label{DistributedEdgeProblemUpdate2}\\
\lambda_{(i,j),i}^{k+1}=&(1-\rho\beta)\lambda_{(i,j),i}^{k}-\rho(A_{ij}x_i^{k+1}+A_{ji}x_j^{k+1}+z_{(i,j),i}^{k+1}), \notag\\
&\forall j\in\mathcal{N}_i, \label{DistributedEdgeProblemUpdate3}
\end{align}
\end{subequations}
where $z_{(i,j),i}$ is a local copy of $z_{(i,j)}$.
This implementation is fully decentralized, as computing the coupling terms only requires data exchange between $i$ and its immediate neighbors.

\begin{figure}[ht]
	\centering
     \vspace{-5pt}
 \subfloat{\includegraphics[width=0.85\linewidth]{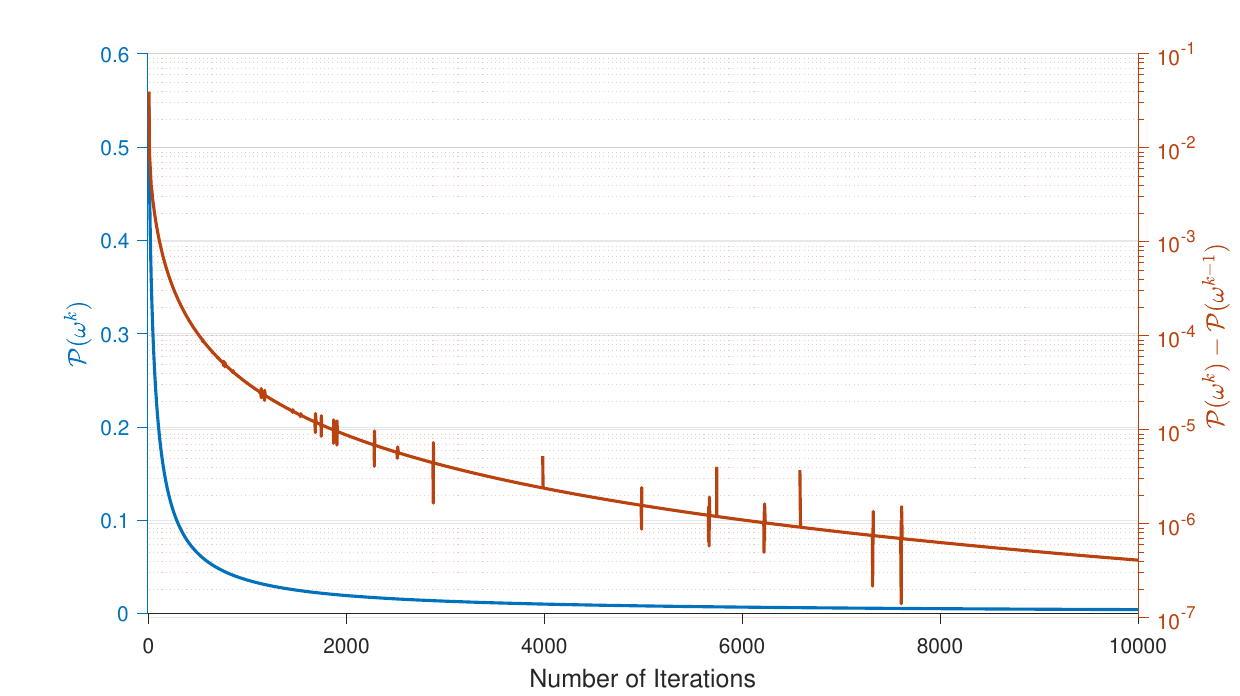}} \\
  \vspace{-5pt}
 \subfloat{\includegraphics[width=0.85\linewidth]{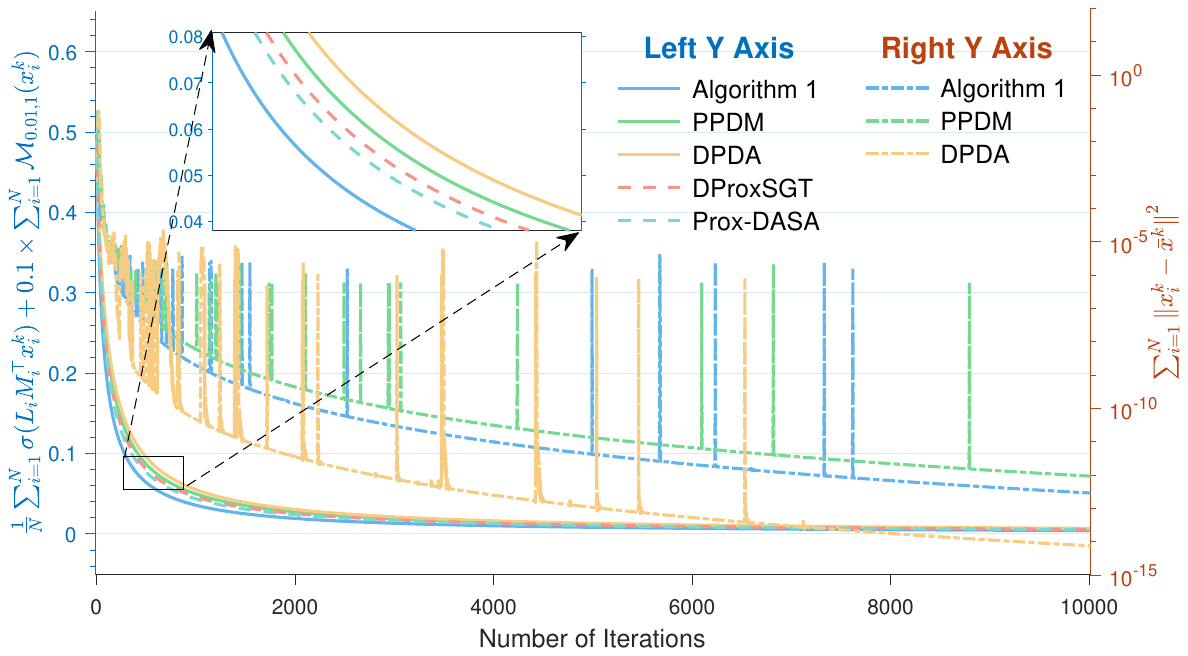}}
 \vspace{-5pt}
	\caption{The trajectory of Lyapunov function, as well as the iteration results of objective function and consensus residuals. 
    Note that the spikes in the trajectories are attributed to the nonsmooth thresholding nature of the SCAD regularizer.}\label{fig_e2}
 \vspace{-15pt}
\end{figure}

We validate the proposed decentralized algorithm \eqref{DistributedEdgeProblemUpdate} on a binary classification task using the MNIST dataset (digits `0' and `8'). 
The optimization problem involves a smooth nonconvex Sigmoid loss function and a nonsmooth weakly convex SCAD regularizer, distributed over a ring network of $N=5$ agents.
For matrices in the partial consensus constraint, they are set to be close to $\mathbf{I}$ and $-\mathbf{I}$.
The experimental results are presented in Fig. \ref{fig_e2}. 
The top panel illustrates the monotonic decrease of the Lyapunov function $\mathcal{P}$, which is consistent with the sufficient descent property established in the theoretical analysis.
The bottom panel compares the convergence of objective values and consensus residuals against state-of-the-art decentralized methods, including PPDM \cite{wang2021distributed}, DPDA \cite{chen2021communication}, DProxSGT \cite{yan2023compressed}, and Prox-DASA \cite{xiao2023one}. 
PPG-ADMM (with $\beta=10^{-8}$) exhibits competitive convergence rate and high solution accuracy, demonstrating its effectiveness in solving practical NCOPs without centralized coordination.

\section{Conclusion}
This paper proposes the PPG-ADMM framework for solving NCOPs. 
By removing the conventional smoothness assumption on the objective function and relaxing the commonly imposed range inclusion condition on constraint matrices, the framework significantly broadens the applicability of ADMM.
The proposed algorithm achieves convergence to an $\epsilon$-approximate stationary point with a sublinear convergence rate of $\mathcal{O}(1/\sqrt{K})$ and an iteration complexity of $\mathcal{O}(1/\epsilon^2)$. 
By selecting a sufficiently small perturbation parameter, it can attain an $\epsilon$-stationary point, ensuring stronger optimality guarantees.
The framework is further extended to solve a practical decentralized optimization problem.
The future research will explore its applications in large-scale networked systems.

\bibliographystyle{IEEEtran}
\bibliography{ref_tn}

\end{document}